\documentclass[11pt]{article}
	\newcommand{\mainTitle}{A formula for systems of Boolean polynomial equations and applications to computational complexity}	
	\newcommand{\authorName}{Tomoya Machide}	
	\newcommand{\organizationNameFst}{Global Research Center for Big Data Mathematics, National Institute of Informatics}
	\newcommand{\placeAddressFst}{2-1-2 Hitotsubashi, Chiyoda-ku, Tokyo 101-8430, Japan}
	\newcommand{\emailAddressFst}{machide@nii.ac.jp}
	\newcommand{\MSCname}{03D15, 03G05, 08A40 (Primary); 06E30, 13P15, 68W30 (Secondary)} 
	\newcommand{\keyWord}{binary tree, Boolean polynomial, computational complexity, formula, system of polynomial equations} 

\usepackage{graphicx}
\usepackage{geometry}               
\usepackage{ascmac}
\usepackage{amsmath}
\usepackage{amssymb}
\usepackage{amsthm}
\usepackage[colorlinks=true]{hyperref}		
\usepackage{arydshln}
\usepackage{multicol}
\usepackage{tikz}

	\DeclareMathOperator*{\OPlus}{\bigoplus}
	\DeclareMathOperator*{\LAnd}{\land}	\DeclareMathOperator*{\LAND}{\bigwedge}	
	\DeclareMathOperator*{\LOr}{\lor}		\DeclareMathOperator*{\LOR}{\bigvee}		
	\DeclareMathOperator*{\dmoPR}{\Pr}
	\newcommand{\nbk}[3]{#1#3#2}		
	\newcommand{\bgbk}[3]{\bigl{#1}#3\bigr{#2}}	
	\newcommand{\Bgbk}[3]{\Bigl{#1}#3\Bigr{#2}}			
	\newcommand{\bggbk}[3]{\biggl{#1}#3\biggr{#2}}			
	\newcommand{\Bggbk}[3]{\Biggl{#1}#3\Biggr{#2}}
	\newcommand{\autobk}[3]{\left#1#3\right#2}
	\newcommand{\nbkD}[5]{#1#2#5#3#4}		
	\newcommand{\bgbkD}[5]{\bigl{#1}\bigl{#2}#5\bigr{#3}\bigr{#4}}	
	\newcommand{\BgbkD}[5]{\Bigl{#1}\Bigl{#2}#5\Bigr{#3}\Bigr{#4}}	
	\newcommand{\bggbkD}[5]{\biggl{#1}\biggl{#2}#5\biggr{#3}\biggr{#4}}	
	\newcommand{\BggbkD}[5]{\Biggl{#1}\Biggl{#2}#5\Biggr{#3}\Biggr{#4}}	
	\newcommand{\autobkD}[5]{\left#1\left#2#5\right#3\right#4}	
	\newcommand{\mcbk}[4][?]{\ifx n#1\nbk{#2}{#3}{#4}\else\ifx b#1\bgbk{#2}{#3}{#4}\else\ifx B#1\Bgbk{#2}{#3}{#4}\else\ifx g#1\bggbk{#2}{#3}{#4}\else\ifx G#1\Bggbk{#2}{#3}{#4}\else\ifx a#1\autobk{#2}{#3}{#4}\else\ifx !#1{#4}\else#4\fi\fi\fi\fi\fi\fi\fi}
	\newcommand{\mcbkD}[4][?]{\ifx n#1\nbkD{#2}{#2}{#3}{#3}{#4}\else\ifx b#1\bgbkD{#2}{#2}{#3}{#3}{#4}\else\ifx B#1\BgbkD{#2}{#2}{#3}{#3}{#4}\else\ifx g#1\bggbkD{#2}{#2}{#3}{#3}{#4}\else\ifx G#1\BggbkD{#2}{#2}{#3}{#3}{#4}\else\ifx a#1\autobkD{#2}{#2}{#3}{#3}{#4}\else\ifx !#1{#4}\else#4\fi\fi\fi\fi\fi\fi\fi}
	\newcommand{\nsgsb}[1]{#1}		
	\newcommand{\bgsgsb}[1]{\big{#1}}	
	\newcommand{\Bgsgsb}[1]{\Big{#1}}			
	\newcommand{\bggsgsb}[1]{\bigg{#1}}			
	\newcommand{\Bggsgsb}[1]{\Bigg{#1}}
	\newcommand{\mcsgsb}[2][?]{\ifx n#1\nsgsb{#2}\else\ifx b#1\bgsgsb{#2}\else\ifx B#1\Bgsgsb{#2}\else\ifx g#1\bggsgsb{#2}\else\ifx G#1\Bggsgsb{#2}\else#2\fi\fi\fi\fi\fi}
	\newcommand{\myEqSpace}{\,} 	\newlength{\myEqSpaceLen} 	
	\newcommand{\eLt}[1]{\widehat{#1}}
	\newcommand{\mLt}[1]{\widetilde{#1}}

	\newcommand{\bkR}[2][n]{\mcbk[#1]{(}{)}{#2}}						
	\newcommand{\bkS}[2][n]{\mcbk[#1]{[}{]}{#2}}						
	\newcommand{\bkB}[2][n]{\mcbk[#1]{\{}{\}}{#2}}						
	\newcommand{\bkAll}[4][n]{\mcbk[#1]{#2}{#3}{#4}}
		
	\newcommand{\nFc}[3][n]{#2\bkR[#1]{#3}}					
				
	\newcommand{\idFc}[4][n]{\id{#2}{#3}\bkR[#1]{#4}}			
	\newcommand{\pwFc}[4][n]{\pw{#2}{#3}\bkR[#1]{#4}}			
			
		\newcommand{\Fc}{\nFc}		
			
	\newcommand{\alp}{\alpha}

	\newcommand{\gam}{\gamma} 
	\newcommand{\Gam}{\Gamma} 
	
	\newcommand{\ep}{\varepsilon}

	\newcommand{\ome}{\omega}

	\newcommand{\bfAlp}[1][s]{\ifx s#1{\boldsymbol\alpha}\else{\boldsymbol??none??}\fi}		
	\newcommand{\bfBeta}[1][s]{\ifx s#1{\boldsymbol\beta}\else{\boldsymbol??none??}\fi}		
	\newcommand{\bfDelta}[1][s]{\ifx s#1{\boldsymbol\delta}\else{\boldsymbol??none??}\fi}		
	\newcommand{\bfGam}[1][s]{\ifx s#1{\boldsymbol\gam}\else{\boldsymbol\Gam}\fi}			
	\newcommand{\bfIota}[1][s]{\ifx s#1{\boldsymbol\iota}\else{\boldsymbol??none??}\fi}			
	\newcommand{\bfA}[1][s]{\ifx s#1{\bf a}\else{\bf A}\fi}
	\newcommand{\bfB}[1][s]{\ifx s#1{\bf b}\else{\bf B}\fi}
	\newcommand{\bfC}[1][s]{\ifx s#1{\bf c}\else{\bf C}\fi}
	\newcommand{\bfD}[1][s]{\ifx s#1{\bf d}\else{\bf D}\fi}
	\newcommand{\bfE}[1][s]{\ifx s#1{\bf e}\else{\bf E}\fi}
	\newcommand{\bfH}[1][s]{\ifx s#1{\bf h}\else{\bf H}\fi}
	\newcommand{\bfI}[1][s]{\ifx s#1{\bf i}\else{\bf I}\fi}
	\newcommand{\bfJ}[1][s]{\ifx s#1{\bf j}\else{\bf J}\fi}
	\newcommand{\bfK}[1][s]{\ifx s#1{\bf k}\else{\bf K}\fi}
	\newcommand{\bfL}[1][s]{\ifx s#1{\bf l}\else{\bf L}\fi}
	\newcommand{\bfM}[1][s]{\ifx s#1{\bf m}\else{\bf M}\fi}
	\newcommand{\bfN}[1][s]{\ifx s#1{\bf n}\else{\bf N}\fi}
	\newcommand{\bfP}[1][s]{\ifx s#1{\bf p}\else{\bf P}\fi}
	\newcommand{\bfQ}[1][s]{\ifx s#1{\bf q}\else{\bf Q}\fi}
	\newcommand{\bfR}[1][s]{\ifx s#1{\bf r}\else{\bf R}\fi}
	\newcommand{\bfS}[1][s]{\ifx s#1{\bf s}\else{\bf S}\fi}
	\newcommand{\bfT}[1][s]{\ifx s#1{\bf t}\else{\bf T}\fi}
	\newcommand{\bfU}[1][s]{\ifx s#1{\bf u}\else{\bf U}\fi}
	\newcommand{\bfV}[1][s]{\ifx s#1{\bf v}\else{\bf V}\fi}
	\newcommand{\bfW}[1][s]{\ifx s#1{\bf w}\else{\bf W}\fi}
	\newcommand{\bfX}[1][s]{\ifx s#1{\bf x}\else{\bf X}\fi}		
	\newcommand{\bfY}[1][s]{\ifx s#1{\bf y}\else{\bf Y}\fi}		
	\newcommand{\bfZ}[1][s]{\ifx s#1{\bf z}\else{\bf Z}\fi}		
	\newcommand{\bfZero}[1][?]{{\bf 0}}
	\newcommand{\bfOne}[1][?]{{\bf 1}}
	\newcommand{\bfEp}[1][s]{\ifx s#1{\boldsymbol\ep}\else{\boldsymbol{\mathcal{E}}}\fi}
	
	\newcommand{\mVert}[1][n]{{\,\mcsgsb[#1]{\vert}\,}}		
	
	\newcommand{\SetO}[2][n]{\bkB[#1]{#2}}
	\newcommand{\SetT}[3][n]{\bkB[#1]{#2\mVert#3}}
		\newcommand{\Set}{\SetO}

	\newcommand{\setZ}{\mathbb{Z}}

	\newcommand{\setF}[1][?]{\ifx ?#1\mathbb{F}\else\mathbb{F}_{#1}\fi}
	
	\newcommand{\setB}{\mathbb{B}}	
	\newcommand{\matI}[1][?]{\ifx #1?I\else I_{#1}\fi}	
	
	\newcommand{\gpKleinF}[1][?]{V}
		
	\newcommand{\gpu}[1][?]{\ifx?#1e\else e_{#1}\fi}		

	\newcommand{\vPack}[1][10]{\vspace{-#1pt}}
	
	\newcommand{\lnA}[1][]{&  &}
	\newcommand{\lnP}[2][1]{\ifx1#1\myEqSpace#2\myEqSpace\else\myEqSpace\myEqSpace#2\myEqSpace\myEqSpace\fi}
	\newcommand{\lnAP}[2][]{& #2 &}
	\newcommand{\lnAH}[1][\nonumber]{#1\\ & &}
	\newcommand{\lnAHP}[2][\nonumber]{#1 \\ & #2 &}

		\newcommand{\slnAH}[1][?]{\\}
		
	\newcommand{\refEq}[1]{(\ref{#1})}	
		
	\newcommand{\pcstSpForRefThm}{\;}		
	\newcommand{\refHL}[2]{#1\pcstSpForRefThm\ref{#2}}		
	\newcommand{\refHLm}[3][?]{\ifx?#1#2\pcstSpForRefThm#3\else#2#3\fi}
	\newcommand{\refThm}[2][?]{\ifx?#1\refHL{Theorem}{#2}\else\ifx s#1\refHL{Theorems}{#2}\else{[argument error]}\fi\fi}
	\newcommand{\refProp}[2][?]{\ifx?#1\refHL{Proposition}{#2}\else\ifx s#1\refHL{Propositions}{#2}\else{[argument error]}\fi\fi}
	\newcommand{\refLem}[2][?]{\ifx?#1\refHL{Lemma}{#2}\else\ifx s#1\refHL{Lemmas}{#2}\else{[argument error]}\fi\fi}
	\newcommand{\refCor}[2][?]{\ifx?#1\refHL{Corollary}{#2}\else\ifx s#1\refHL{Corollaries}{#2}\else{[argument error]}\fi\fi}
	\newcommand{\refDef}[2][?]{\ifx?#1\refHL{Definition}{#2}\else\ifx s#1\refHL{Definitions}{#2}\else{[argument error]}\fi\fi}
	\newcommand{\refRem}[2][?]{\ifx?#1\refHL{Remark}{#2}\else\ifx s#1\refHL{Remarks}{#2}\else{[argument error]}\fi\fi}
	\newcommand{\refTab}[2][?]{\ifx?#1\refHL{Table}{#2}\else\ifx s#1\refHL{Tables}{#2}\else{[argument error]}\fi\fi}
	\newcommand{\refFig}[2][?]{\ifx?#1\refHL{Figure}{#2}\else\ifx s#1\refHL{Figures}{#2}\else{[argument error]}\fi\fi}
	\newcommand{\refSec}[2][?]{\ifx?#1\refHL{Section}{#2}\else\ifx s#1\refHL{Sections}{#2}\else{[argument error]}\fi\fi}
	\newcommand{\refApp}[2][?]{\ifx?#1\refHL{Appendix}{#2}\else\ifx s#1\refHL{Appendices}{#2}\else{[argument error]}\fi\fi}
	\newcommand{\refAlg}[2][?]{\ifx?#1\refHL{Algorithm}{#2}\else\ifx s#1\refHL{Algorithms}{#2}\else{[argument error]}\fi\fi}
	\newcommand{\refPrc}[2][?]{\ifx?#1\refHL{Process}{#2}\else\ifx s#1\refHL{Processes}{#2}\else{[argument error]}\fi\fi}
		\newcommand{\refSect}{\refSec}
		\newcommand{\refPrp}{\refProp}
	\newcommand{\refThmA}[2][?]{\ifx?#1\refHLm{Theorem}{#2}\else \refHLm{Theorems}{#2}\fi}
	\newcommand{\refPropA}[2][?]{\ifx?#1 \refHLm[#1]{Proposition}{#2}\else \refHLm{Propositions}{#2}\fi}
	\newcommand{\refLemA}[2][?]{\ifx?#1\refHLm{Lemma}{#2}\else \refHLm{Lemmas}{#2}\fi}
	\newcommand{\refCorA}[2][?]{\ifx?#1\refHLm{Corollary}{#2}\else \refHLm{Corollaries}{#2}\fi}
	\newcommand{\refDefA}[2][?]{\ifx?#1\refHLm{Definition}{#2}\else \refHLm{Definitionss}{#2}\fi}
	\newcommand{\refRemA}[2][?]{\ifx?#1\refHLm{Remark}{#2}\else \refHLm{Remarks}{#2}\fi}
	\newcommand{\refSecA}[2][?]{\ifx?#1\refHLm{Section}{#2}\else \refHLm{Sections}{#2}\fi}
	\newcommand{\refAppA}[2][?]{\ifx?#1\refHLm{Appendix}{#2}\else \refHLm{Appendices}{#2}\fi}
	\newcommand{\refAlgA}[2][?]{\ifx?#1\refHLm{Algorithm}{#2}\else \refHLm{Algorithms}{#2}\fi}
	\newcommand{\refPrcA}[2][?]{\ifx?#1\refHLm{Process}{#2}\else \refHLm{Processes}{#2}\fi}

	\newcommand{\frc}[3][?]{\ifx s#1#3/#2\else\ifx b#1(#3)/#2\else\ifx d#1\dfrac{#3}{#2}\else\ifx t#1\tfrac{#3}{#2}\else\frac{#3}{#2}\fi\fi\fi\fi}
	\newcommand{\bnm}[3][?]{\binom{#3}{#2}}

	\newcommand{\pw}[3][?]{\ifx!#3{#2}^{#3}\else#2^{#3}\fi}
	\newcommand{\id}[3][?]{#2_{#3}}
	
	\newcommand{\pwR}[3][a]{\ifx!#1{\bkR[#1]{#2}}^{#3}\else\bkR[#1]{#2}^{#3}\fi}
	\newcommand{\pwB}[3][a]{\ifx!#1{\bkB[#1]{#2}}^{#3}\else\bkB[#1]{#2}^{#3}\fi}
	\newcommand{\pwS}[3][a]{\ifx!#1{\bkS[#1]{#2}}^{#3}\else\bkS[#1]{#2}^{#3}\fi}

	\newcommand{\nIntO}[2][?]{\ifx l#1\int\limits_{#2}\else\ifx t#1{\textstyle\int\limits_{#2}}\else\int_{#2}\fi\fi}
	\newcommand{\nIntT}[3][?]{\ifx l#1\int\limits_{#2}^{#3}\else\if t#1{\textstyle\int\limits_{#2}^{#3}}\else\int_{#2}^{#3}\fi\fi}	
	\newcommand{\nIntN}[1][?]{\ifx l#1\int\limits\else\ifx t#1{\textstyle\int\limits}\else\int\fi\fi}

	\newcommand{\tpT}[3][a]{ {#2}\atop \bkR[#1]{#3} }

	\newcommand{\nSmO}[2][?]{\ifx l#1\sum\limits_{#2}\else\ifx t#1{\textstyle\sum\limits_{#2}}\else\sum_{#2}\fi\fi}
	\newcommand{\nSmT}[3][?]{\ifx l#1\sum\limits_{#2}^{#3}\else\if t#1{\textstyle\sum\limits_{#2}^{#3}}\else\sum_{#2}^{#3}\fi\fi}	
	\newcommand{\nSmN}[1][?]{\ifx l#1\sum\limits\else\ifx t#1{\textstyle\sum\limits}\else\sum\fi\fi}
	\newcommand{\pSm}[2][?]{\ifx t#1 \sum_{#2}^{\prime} \else \sideset{}{^\prime}\sum_{#2} \fi}
	\newcommand{\pSmT}[3][?]{\ifx t#1 \sum_{#2}^{\prime#3} \else \sideset{}{^\prime}\sum_{#2}^{#3} \fi}	
	\newcommand{\pSmN}[1][?]{\ifx t#1 \sum^{\prime} \else \sideset{}{^\prime}\sum \fi}
	\newcommand{\dSm}[2][?]{\ifx t#1 \sum_{#2}^{\dagger} \else \sideset{}{^\dagger}\sum_{#2} \fi}
	\newcommand{\dSmT}[3][?]{\ifx t#1 \sum_{#2}^{\dagger#3} \else \sideset{}{^\dagger}\sum_{#2}^{#3} \fi}	
	\newcommand{\dSmN}[1][?]{\ifx t#1 \sum^{\dagger} \else \sideset{}{^\dagger}\sum \fi}
	\newcommand{\tpTSm}[3][?]{\nSmO[#1]{\tpT{#2}{#3}}}

		\newcommand{\Sm}{\nSmO}			\newcommand{\SmT}{\nSmT}			
		\newcommand{\tpSm}{\tpTSm}
		
	\newcommand{\nPd}[2][?]{\ifx l#1 \prod\limits_{#2} \else \prod_{#2} \fi}
	\newcommand{\nPdT}[3][?]{\ifx l#1 \prod\limits_{#2}^{#3} \else \prod_{#2}^{#3} \fi}	
	
	\newcommand{\tpTPdT}[4][?]{ \nPdT[#1]{ \tpT{#2}{#3} }{#4} }
		\newcommand{\Pd}{\nPd}
		\newcommand{\PdT}{\nPdT}
		
		\newcommand{\tpPdT}{\tpTPdT}
	\newcommand{\nOPs}[2][?]{\ifx l#1 \OPlus\limits_{#2} \else \OPlus_{#2} \fi}
	\newcommand{\nOPsT}[3][?]{\ifx l#1 \OPlus\limits_{#2}^{#3} \else \OPlus_{#2}^{#3} \fi}	
	\newcommand{\pOPs}[2][?]{\ifx t#1 \OPlus_{#2}^{\prime} \else \sideset{}{^\prime}\OPlus_{#2} \fi}
	\newcommand{\pOPsT}[3][?]{\ifx t#1 \OPlus_{#2}^{\prime#3} \else \sideset{}{^\prime}\OPlus_{#2}^{#3} \fi}

	\newcommand{\nAD}[2][?]{\ifx l#1 \LAND\limits_{#2} \else \LAND_{#2} \fi}
	\newcommand{\nADT}[3][?]{\ifx l#1 \LAND\limits_{#2}^{#3} \else \LAND_{#2}^{#3} \fi}	
	\newcommand{\pAD}[2][?]{\ifx t#1 \LAND_{#2}^{\prime} \else \sideset{}{^\prime}\LAND_{#2} \fi}
	\newcommand{\pADT}[3][?]{\ifx t#1 \LAND_{#2}^{\prime#3} \else \sideset{}{^\prime}\LAND_{#2}^{#3} \fi}

		\newcommand{\AD}{\nAD}

	\newcommand{\nAd}[2][?]{\ifx l#1 \LAnd\limits_{#2} \else \LAnd_{#2} \fi}
	\newcommand{\nAdT}[3][?]{\ifx l#1 \LAnd\limits_{#2}^{#3} \else \LAnd_{#2}^{#3} \fi}	
	\newcommand{\pAd}[2][?]{\ifx t#1 \LAnd_{#2}^{\prime} \else \sideset{}{^\prime}\LAnd_{#2} \fi}
	\newcommand{\pAdT}[3][?]{\ifx t#1 \LAnd_{#2}^{\prime#3} \else \sideset{}{^\prime}\LAnd_{#2}^{#3} \fi}

		\newcommand{\Ad}{\nAd}

	\newcommand{\nOR}[2][?]{\ifx l#1 \LOR\limits_{#2} \else \LOR_{#2} \fi}
	\newcommand{\nORT}[3][?]{\ifx l#1 \LOR\limits_{#2}^{#3} \else \LOR_{#2}^{#3} \fi}	
	\newcommand{\pOR}[2][?]{\ifx t#1 \LOR_{#2}^{\prime} \else \sideset{}{^\prime}\LOR_{#2} \fi}
	\newcommand{\pORT}[3][?]{\ifx t#1 \LOR_{#2}^{\prime#3} \else \sideset{}{^\prime}\LOR_{#2}^{#3} \fi}

		\newcommand{\OR}{\nOR}

	\newcommand{\nOr}[2][?]{\ifx l#1 \LOr\limits_{#2} \else \LOr_{#2} \fi}
	\newcommand{\nOrT}[3][?]{\ifx l#1 \LOr\limits_{#2}^{#3} \else \LOr_{#2}^{#3} \fi}	
	\newcommand{\pOr}[2][?]{\ifx t#1 \LOr_{#2}^{\prime} \else \sideset{}{^\prime}\LOr_{#2} \fi}
	\newcommand{\pOrT}[3][?]{\ifx t#1 \LOr_{#2}^{\prime#3} \else \sideset{}{^\prime}\LOr_{#2}^{#3} \fi}

		\newcommand{\Or}{\nOr}

	\newcommand{\nPR}[2][?]{\ifx l#1 \dmoPR\limits_{#2} \else \dmoPR_{#2} \fi}
	\newcommand{\nPRT}[3][?]{\ifx l#1 \dmoPR\limits_{#2}^{#3} \else \dmoPR_{#2}^{#3} \fi}	
	\newcommand{\pPR}[2][?]{\ifx t#1 \dmoPR_{#2}^{\prime} \else \sideset{}{^\prime}\dmoPR_{#2} \fi}
	\newcommand{\pPRT}[3][?]{\ifx t#1 \dmoPR_{#2}^{\prime#3} \else \sideset{}{^\prime}\dmoPR_{#2}^{#3} \fi}

	\newcommand{\nIs}[2][?]{\ifx l#1 \bigcap\limits_{#2}\else\ifx b#1 \bigcap_{#2}\else{\textstyle\bigcap\limits_{#2}}\fi\fi}
	\newcommand{\nIsT}[3][?]{\ifx l#1 \bigcap\limits_{#2}^{#3}\else\ifx b#1 \bigcap_{#2}^{#3}\else{\textstyle\bigcap\limits_{#2}^{#3}}\fi\fi}	
	\newcommand{\pIs}[2][?]{\ifx t#1 \bigcap_{#2}^{\prime} \else \sideset{}{^\prime}\bigcap_{#2} \fi}
	\newcommand{\pIsT}[3][?]{\ifx t#1 \bigcap_{#2}^{\prime#3} \else \sideset{}{^\prime}\bigcap_{#2}^{#3} \fi}

	\newcommand{\nUn}[2][?]{\ifx L#1 \bigcup\limits_{#2}\else\ifx b#1 \bigcup_{#2}\else\ifx l#1{\textstyle\bigcup\limits_{#2}}\else{\textstyle\bigcup_{#2}}\fi\fi\fi}
	\newcommand{\nUnT}[3][?]{\ifx L#1 \bigcup\limits_{#2}^{#3}\else\ifx b#1 \bigcup_{#2}^{#3}\else\ifx l#1{\textstyle\bigcup\limits_{#2}^{#3}}\else{\textstyle\bigcup_{#2}^{#3}}\fi\fi\fi}	
	\newcommand{\pUn}[2][?]{\ifx t#1 \bigcup_{#2}^{\prime} \else \sideset{}{^\prime}\bigcup_{#2} \fi}
	\newcommand{\pUnT}[3][?]{\ifx t#1 \bigcup_{#2}^{\prime#3} \else \sideset{}{^\prime}\bigcup_{#2}^{#3} \fi}	
	\newcommand{\dUn}[2][?]{\ifx L#1 \bigsqcup\limits_{#2}\else\ifx b#1 \bigsqcup_{#2}\else\ifx l#1{\textstyle\bigsqcup\limits_{#2}}\else{\textstyle\bigsqcup{#2}}\fi\fi\fi}
	\newcommand{\dUnT}[3][?]{\ifx L#1 \bigsqcup\limits_{#2}^{#3}\else\ifx b#1 \bigsqcup_{#2}^{#3}\else\ifx l#1{\textstyle\bigsqcup\limits_{#2}^{#3}}\else{\textstyle\bigsqcup{#2}^{#3}}\fi\fi\fi}

	\newcommand{\nLm}[2][?]{\ifx l#1 \lim\limits_{#2} \else \lim_{#2} \fi}
	\newcommand{\iLm}[2][?]{\ifx l#1 \liminf\limits_{#2} \else \liminf_{#2} \fi}
	\newcommand{\sLm}[2][?]{\ifx l#1 \limsup\limits_{#2} \else \limsup_{#2} \fi}
		
	\newcommand{\ntpMaxT}[3][?]{ \underset{#2}{\max\;}\bkB[#1]{#3} }

		\newcommand{\tpMaxT}{\ntpMaxT}		\newcommand{\tpMax}{\tpMaxT}

	\newcommand{\glcondEnvLineHead}[1]{ \ifx*#1 \begin{eqnarray*} \else \begin{eqnarray}  \label{#1} \fi }
	\newcommand{\glcondEnvLineTail}[1]{ \ifx*#1 \end{eqnarray*} \else \end{eqnarray} \fi }
	\newcommand{\glcondDis}[1]{\ifx d#1 \displaystyle \fi}
	\newcommand{\glcmdEqShift}{\hspace{-20pt}}
	\newcommand{\glcmdHLineCWiden}{\rule{0cm}{15pt}}	\newcommand{\glcdH}{\glcmdHLineCWiden}
	\newcommand{\lccondPar}[1]{\ifx#1p \\ \fi}

		\newcommand{\envMO}[2][*]{$\ifx d#1 \displaystyle \fi#2$}
		\newcommand{\envMT}[3][*]{$\ifx d#1 \displaystyle \fi#2=#3$}
		\newcommand{\envMTDef}[3][*]{$\ifx d#1 \displaystyle \fi#2:=#3$}
		\newcommand{\envMTPt}[4][*]{$\ifx d#1 \displaystyle \fi#3#2#4$}
			\newcommand{\envM}{\envMT}
			
			\newcommand{\envMPt}{\envMTPt}
		\newcommand{\envMTh}[4][*]{$\ifx d#1 \displaystyle \fi#2=#3=#4$}
		\newcommand{\envMThDef}[4][*]{$\ifx d#1 \displaystyle \fi#2:=#3=#4$}
		\newcommand{\envMThPt}[5][*]{$\ifx d#1 \displaystyle \fi#3#2#4#2#5$}
		\newcommand{\envMThPte}[6][*]{$\ifx d#1 \displaystyle \fi#2#3#4#5#6$}
		\newcommand{\envMF}[5][*]{$\ifx d#1 \displaystyle \fi#2=#3=#4=#5$}
		\newcommand{\envMFPt}[6][*]{$\ifx d#1 \displaystyle \fi#3#2#4#2#5#2#6$}
		
		\newcommand{\envLineT}[3][*]{ \glcondEnvLineHead{#1} & &\glcmdEqShift#2\nonumber\\&=&#3 \glcondEnvLineTail{#1} }

			\newcommand{\envLine}{\envLineT}

		\newcommand{\envHLineT}[3][*]{ \glcondEnvLineHead{#1} #2&=&#3\glcondEnvLineTail{#1} }
		\newcommand{\envHLineTDef}[3][*]{ \glcondEnvLineHead{#1} #2&:=&#3\glcondEnvLineTail{#1} }
		\newcommand{\envHLineTPt}[4][*]{\glcondEnvLineHead{#1} #3&#2&#4\glcondEnvLineTail{#1}}
		
			\newcommand{\envHLine}{\envHLineT}
			\newcommand{\envHLineDef}{\envHLineTDef}
			\newcommand{\envHLinePt}{\envHLineTPt}
			
		\newcommand{\envHLineF}[5][*]{ \glcondEnvLineHead{#1} #2&=&#3\nonumber\\&=&#4\nonumber\\&=&#5 \glcondEnvLineTail{#1}}

		\newcommand{\envPLine}[2][*]{\glcondEnvLineHead{#1} #2\glcondEnvLineTail{#1}}
		
	\newcommand{\lcparaCase}{\vspace{3pt}}

	\newcommand{\matu}[1][?]{\ifx#1?I\else I_{#1}\fi}

	\newcommand{\vlA}[2][n]{\bkAll[#1]{|}{|}{#2}}
				\newcommand{\abs}{\vlA}
		
	\newcommand{\nmSet}[2][n]{\vlA[#1]{#2}}	
		
	\newcommand{\sbLn}[2][n]{\Fc[#1]{O}{#2}}		
	\newcommand{\sbLs}[2][n]{\Fc[#1]{o}{#2}}		
	\newcommand{\sbLp}[2][n]{\pwFc[#1]{O}{\star}{#2}}	

		\newcommand{\sbL}{\sbLn}

	\newcommand{\sTx}[2][?]{ \ifx t#1{\tiny #2} \else \ifx s#1{\scriptsize #2} \else \ifx f#1{\footnotesize #2} \else \ifx S#1{\small #2} \else \ifx n#1{\normalsize #2} \else \ifx l#1{\large #2} \else \ifx L#1{\Large #2} \else \ifx R#1{\LARGE #2} \else \ifx h#1{\huge #2} \else \ifx H#1{\Huge #2} \else \ifx ?#1 #2 \else #2 \fi\fi\fi\fi\fi\fi\fi\fi\fi\fi\fi }

	\newcommand{\osMTx}[3][?]{\overset{#3}{#2}}
	\newcommand{\usbMTx}[3][?]{\underset{#2}{\underbrace{#3}}}

		\newcommand{\osTx}{\osMTx}
		\newcommand{\usbTx}{\usbMTx}
		
	\newcommand{\raTx}[3][?]{\raisebox{#2pt}[0pt][0pt]{\ifx d#1\displaystyle\fi#3}}
	\newcommand{\raMTx}[3][?]{\raisebox{#2pt}[0pt][0pt]{$\ifx d#1\displaystyle\fi#3$}}
	
	\newcommand{\roTx}[3][?]{\rotatebox[origin=c]{#2}{#3}}

	\newcommand{\envCenter}[2][*]{\ifx*#1\begin{center}\else\begin{center}[#1]\fi #2\end{center}}
	\newcommand{\envFlushleft}[2][*]{\ifx*#1\begin{flushleft}\else\begin{flushleft}[#1]\fi #2\end{flushleft}}
	\newcommand{\envFlushright}[2][*]{\ifx*#1\begin{flushright}\else\begin{flushright}[#1]\fi #2\end{flushright}}
		
	\newcommand{\envItemIm}[2][*]{\ifx*#1\begin{itemize}\else\begin{itemize}[#1]\fi #2\end{itemize}}
	\newcommand{\envItemDp}[2][*]{\ifx*#1\begin{description}\else\begin{description}[#1]\fi #2\end{description}}
	\newcommand{\envItemEm}[2][*]{\ifx*#1\begin{enumerate}\else\begin{enfumerate}[#1]\fi #2\end{enumerate}}
		\newcommand{\envItem}{\envItemIm}
		
	\newcommand{\envMultCol}[3][*]{\ifx1#2#3\else\begin{multicols}{#2}\ifx*#1\else\mbox{}\vspace{-#1pt}\fi#3\end{multicols}\fi}
	
	\newcommand{\envPic}[2][*]{\begin{picture}#2\end{picture}}

\theoremstyle{plain}
\newtheorem{theorem}{THEOREM}[section]
\newtheorem{proposition}[theorem]{PROPOSITION}
\newtheorem{lemma}[theorem]{LEMMA}
\newtheorem{corollary}[theorem]{COROLLARY}
\theoremstyle{definition}

\theoremstyle{remark}

\theoremstyle{plain}

\theoremstyle{definition}

\theoremstyle{remark}

\theoremstyle{plain}

\theoremstyle{definition}

\theoremstyle{remark}

\allowdisplaybreaks[4]
\numberwithin{equation}{section}

	\newcommand{\lccondBibitem}[3][]{ \if ?#2 \bibitem{#3} \else \bibitem[#2]{#3} \fi}
	\newcommand{\refPaper}[8][?]{
			\lccondBibitem{#1}{#2}
				#3,			
				\emph{#4}, 	
				#5\ 			
				{\bf #6}		
				(#7),			
				#8.			
		}
	\newcommand{\refPreprint}[6][?]{
			\lccondBibitem{#1}{#2}
				#3,			
				\emph{#4}, 	
				preprint; #5,	
				#6.			
		}

	\newcommand{\refBook}[7][?]{
			\lccondBibitem{#1}{#2}
				#3,			
				\emph{#4}, 	
				#5,			
				#6,			
				#7.			
		}
	\newcommand{\refPaperAlm}[5][?]{
			\lccondBibitem{#1}{#2}
				#3,	 		
				\emph{#4}, 	
				#5		
		}

	\newcommand{\etalTx}[2][?]{#2 \emph{et al.}\!}

	\newcommand{\glcondEnvLineTailPd}[1]{.\ifx*#1 \end{eqnarray*} \else \end{eqnarray} \fi  }
	\newcommand{\glcondEnvLineTailCm}[1]{,\ifx*#1 \end{eqnarray*} \else \end{eqnarray} \fi }
	\newcommand{\prcondEnvEqSpHead}[1]{ \ifx*#1 \begin{equation*}[ERROR] \else \begin{equation}  \label{#1} \fi  }
	\newcommand{\prcondEnvEqSpTail}[1]{\ifx*#1 [ERROR]\end{equation*} \else \end{equation} \fi }

	\newcommand{\envProof}[2][?]{ \par\mbox{}\vspace{-5pt}\\ \ifx?#1\emph{Proof.}\else\emph{Proof of #1.}\fi \ #2 \hfill $\Box$\\ \par}
	
		\newcommand{\envLineTPd}[3][*]{ \glcondEnvLineHead{#1} & &\glcmdEqShift#2\nonumber\\&=&#3 \glcondEnvLineTailPd{#1} }

			\newcommand{\envLinePd}{\envLineTPd}

		\newcommand{\envLineTCmPt}[4][*]{\glcondEnvLineHead{#1} & &\glcmdEqShift#3\nonumber\\&#2&#4 \glcondEnvLineTailCm{#1}}

			\newcommand{\envLineCmPt}{\envLineTCmPt}

		\newcommand{\envLineThPd}[4][*]{ \glcondEnvLineHead{#1} & &\glcmdEqShift#2\nonumber\\&=&#3\nonumber \\&=&#4 \glcondEnvLineTailPd{#1} }
		\newcommand{\envLineThCm}[4][*]{ \glcondEnvLineHead{#1} & &\glcmdEqShift#2\nonumber\\&=&#3\nonumber \\&=&#4 \glcondEnvLineTailCm{#1} }

		\newcommand{\envHLineTPd}[3][*]{ \glcondEnvLineHead{#1} #2&=&#3\glcondEnvLineTailPd{#1} }
		\newcommand{\envHLineTDefPd}[3][*]{ \glcondEnvLineHead{#1} #2&:=&#3\glcondEnvLineTailPd{#1} }
		\newcommand{\envHLineTCm}[3][*]{ \glcondEnvLineHead{#1} #2&=&#3\glcondEnvLineTailCm{#1} }
		\newcommand{\envHLineTCmDef}[3][*]{ \glcondEnvLineHead{#1} #2&:=&#3\glcondEnvLineTailCm{#1} }
		\newcommand{\envHLineTCmPt}[4][*]{\glcondEnvLineHead{#1} #3&#2&#4\glcondEnvLineTailCm{#1}}					
		\newcommand{\envHLineTPdPt}[4][*]{\glcondEnvLineHead{#1} #3&#2&#4\glcondEnvLineTailPd{#1}}

			\newcommand{\envHLinePd}{\envHLineTPd}
			\newcommand{\envHLineDefPd}{\envHLineTDefPd}
			\newcommand{\envHLineCm}{\envHLineTCm}
			\newcommand{\envHLineCmDef}{\envHLineTCmDef}
			\newcommand{\envHLineCmPt}{\envHLineTCmPt}
			\newcommand{\envHLinePdPt}{\envHLineTPdPt}

		\newcommand{\envHLineThPd}[4][*]{ \glcondEnvLineHead{#1} #2&=&#3\nonumber\\&=&#4\glcondEnvLineTailPd{#1}	 }
		
		\newcommand{\envHLineThCm}[4][*]{ \glcondEnvLineHead{#1} #2&=&#3\nonumber\\&=&#4\glcondEnvLineTailCm{#1}}

		\newcommand{\envHLineFPd}[5][*]{ \glcondEnvLineHead{#1} #2&=&#3\nonumber\\&=&#4\nonumber \\&=&#5 \glcondEnvLineTailPd{#1} }

		\newcommand{\envHLineFiCm}[6][*]{ \glcondEnvLineHead{#1} #2&=&#3\nonumber\\&=&#4\nonumber \\&=&#5\nonumber  \\&=&#6\glcondEnvLineTailCm{#1}}

		\newcommand{\envHLineCFCmNme}[5][*]{\begin{eqnarray} #2&=&#3,\\\glcdH#4&=&#5 \glcondEnvLineTailCm{?} }
		\newcommand{\envHLineCFNmePd}[5][*]{\begin{eqnarray} #2&=&#3,\\\glcdH#4&=&#5 \glcondEnvLineTailPd{?} }
		\newcommand{\envHLineCFCmDefNme}[5][*]{\begin{eqnarray} #2&:=&#3,\\\glcdH#4&:=&#5 \glcondEnvLineTailCm{?} }
		\newcommand{\envHLineCFDefNmePd}[5][*]{\begin{eqnarray} #2&:=&#3,\\\glcdH#4&:=&#5 \glcondEnvLineTailPd{?} }

		\newcommand{\envHLineCFNmePdPt}[6][*]{\begin{eqnarray}#3&#2&#4,\\\glcdH#5&#2&#6\glcondEnvLineTailPd{?}}
		\newcommand{\envHLineCFCmNmePt}[6][*]{\begin{eqnarray}#3&#2&#4,\\\glcdH#5&#2&#6\glcondEnvLineTailCm{?}}
		\newcommand{\envHLineCFNmePdPte}[7][*]{\begin{eqnarray}#2&#3&#4,\\\glcdH#5&#6&#7\glcondEnvLineTailPd{?}}
		\newcommand{\envHLineCFCmNmePte}[7][*]{\begin{eqnarray}#2&#3&#4,\\\glcdH#5&#6&#7\glcondEnvLineTailCm{?}}

		\newcommand{\envHLineCSNmePd}[7][*]{\begin{eqnarray} #2&=&#3,\\\glcdH#4&=&#5,\\\glcdH#6&=&#7\glcondEnvLineTailPd{?}}
		\newcommand{\envHLineCSDefNmePd}[7][*]{\begin{eqnarray} #2&:=&#3,\\\glcdH#4&:=&#5,\\\glcdH#6&:=&#7\glcondEnvLineTailPd{?}}
		\newcommand{\envHLineCSCmNme}[7][*]{\begin{eqnarray} #2&=&#3,\\\glcdH#4&=&#5,\\\glcdH#6&=&#7\glcondEnvLineTailCm{?}}
		\newcommand{\envHLineCSCmDefNme}[7][*]{\begin{eqnarray} #2&:=&#3,\\\glcdH#4&:=&#5,\\\glcdH#6&:=&#7\glcondEnvLineTailCm{?}}

		\newcommand{\envHLineCSNmePdPt}[8][*]{\begin{eqnarray}#3&#2&#4,\\\glcdH#5&#2&#6,\\\glcdH#7&#2&#8\glcondEnvLineTailPd{?}}
		\newcommand{\envHLineCSCmNmePt}[8][*]{\begin{eqnarray}#3&#2&#4,\\\glcdH#5&#2&#6,\\\glcdH#7&#2&#8\glcondEnvLineTailCm{?}}
		\newcommand{\envHLineCSNmePdPte}[9][*]{\begin{eqnarray}#2&#3&#4,\\\glcdH#5&#6&#7,\\\glcdH#8&#2&#9\glcondEnvLineTailPd{?}}
		\newcommand{\envHLineCSCmNmePte}[9][*]{\begin{eqnarray}#2&#3&#4,\\\glcdH#5&#6&#7,\\\glcdH#8&#2&#9\glcondEnvLineTailCm{?}}
		
		\newcommand{\envHLineCENmePd}[9][*]{\begin{eqnarray} #2&=&#3,\\\glcdH#4&=&#5,\\\glcdH#6&=&#7,\\\glcdH#8&=&#9\glcondEnvLineTailPd{?}}
		\newcommand{\envHLineCEDefNmePd}[9][*]{\begin{eqnarray} #2&:=&#3,\\\glcdH#4&:=&#5,\\\glcdH#6&:=&#7,\\\glcdH#8&:=&#9\glcondEnvLineTailPd{?}}
		\newcommand{\envHLineCECmNme}[9][*]{\begin{eqnarray} #2&=&#3,\\\glcdH#4&=&#5,\\\glcdH#6&=&#7,\\\glcdH#8&=&#9\glcondEnvLineTailCm{?}}
		\newcommand{\envHLineCECmDefNme}[9][*]{\begin{eqnarray} #2&:=&#3,\\\glcdH#4&:=&#5,\\\glcdH#6&:=&#7,\\\glcdH#8&:=&#9\glcondEnvLineTailCm{?}}

			\newcommand{\pccondPaForPar}[1]{\ifx#1p \\\glcdH \fi}
			\newcommand{\pccondPaForNonnum}[1]{\ifx#1* \nonumber \fi}

		\newcommand{\envPLinePd}[2][*]{\glcondEnvLineHead{#1} #2\glcondEnvLineTailPd{#1}}
		\newcommand{\envPLineCm}[2][*]{\glcondEnvLineHead{#1} #2\glcondEnvLineTailCm{#1}}
	
		\newcommand{\envOTLineCm}[4][*]{\glcondEnvLineHead{#1} #2\lnAP{=}#3\lnP{=}#4,\glcondEnvLineTail{#1}}

			\newcommand{\envOTLineThCm}{\envOTLineCm}

		\newcommand{\envOFLineCm}[5][*]{\glcondEnvLineHead{#1} #2\lnAP{=}#3\lnP{=}#4\lnP{=}#5,\glcondEnvLineTail{#1}}

			\newcommand{\envOFLineFCm}{\envOFLineCm}

		\newcommand{\envMOCm}[2][*]{$\ifx d#1 \displaystyle \fi#2$,}
		\newcommand{\envMOPd}[2][*]{$\ifx d#1 \displaystyle \fi#2$.}
		
		\newcommand{\envMTCm}[3][*]{$\ifx d#1 \displaystyle \fi#2=#3$,}
		\newcommand{\envMTPd}[3][*]{$\ifx d#1 \displaystyle \fi#2=#3$.}
		\newcommand{\envMTCmDef}[3][*]{$\ifx d#1 \displaystyle \fi#2:=#3$,}
		\newcommand{\envMTDefPd}[3][*]{$\ifx d#1 \displaystyle \fi#2:=#3$.}
		\newcommand{\envMTCmPt}[4][*]{$\ifx d#1 \displaystyle \fi#3#2#4$,}
		\newcommand{\envMTPdPt}[4][*]{$\ifx d#1 \displaystyle \fi#3#2#4$.}
			\newcommand{\envMCm}{\envMTCm}
			\newcommand{\envMPd}{\envMTPd}

			\newcommand{\envMCmPt}{\envMTCmPt}
			\newcommand{\envMPdPt}{\envMTPdPt}
		\newcommand{\envMThCm}[4][*]{$\ifx d#1 \displaystyle \fi#2=#3=#4$,}
		\newcommand{\envMThPd}[4][*]{$\ifx d#1 \displaystyle \fi#2=#3=#4$.}
		\newcommand{\envMThCmPt}[5][*]{$\ifx d#1 \displaystyle \fi#3#2#4#2#5$,}
		\newcommand{\envMThPdPt}[5][*]{$\ifx d#1 \displaystyle \fi#3#2#4#2#5$.}
		\newcommand{\envMFCm}[5][*]{$\ifx d#1 \displaystyle \fi#2=#3=#4=#5$,}
		\newcommand{\envMFPd}[5][*]{$\ifx d#1 \displaystyle \fi#2=#3=#4=#5$.}
		\newcommand{\envMFCmPt}[6][*]{$\ifx d#1 \displaystyle \fi#3#2#4#2#5#2#6$,}
		\newcommand{\envMFPdPt}[6][*]{$\ifx d#1 \displaystyle \fi#3#2#4#2#5#2#6$.}

		\newcommand{\prcondHLCPNm}{\hspace{-1pt}}
		\newcommand{\envHLineCFCmNm}[5][*]{ \begin{equation}\begin{split} \ifx*#1 \text{[ERROR;need label name]} \else \label{#1} \fi #2&\prcondHLCPNm\lnP{=}\prcondHLCPNm#3,\\#4&\prcondHLCPNm\lnP{=}\prcondHLCPNm#5, \end{split}\end{equation} }
		\newcommand{\envHLineCFNm}[5][*]{ \begin{equation}\begin{split} \ifx*#1 \text{[ERROR;need label name]} \else \label{#1} \fi #2&\prcondHLCPNm\lnP{=}\prcondHLCPNm#3\\#4&\prcondHLCPNm\lnP{=}\prcondHLCPNm#5, \end{split}\end{equation} }
		\newcommand{\envHLineCFNmPd}[5][*]{ \begin{equation}\begin{split} \ifx*#1 \text{[ERROR;need label name]} \else \label{#1} \fi #2&\prcondHLCPNm\lnP{=}\prcondHLCPNm#3,\\#4&\prcondHLCPNm\lnP{=}\prcondHLCPNm#5. \end{split}\end{equation} }
		\newcommand{\envHLineCFCmDefNm}[5][*]{ \begin{equation}\begin{split} \ifx*#1 \text{[ERROR;need label name]} \else \label{#1} \fi #2&\prcondHLCPNm\lnP{:=}\prcondHLCPNm#3,\\#4&\prcondHLCPNm\lnP{:=}\prcondHLCPNm#5, \end{split}\end{equation} }
		\newcommand{\envHLineCFDefNm}[5][*]{ \begin{equation}\begin{split} \ifx*#1 \text{[ERROR;need label name]} \else \label{#1} \fi #2&\prcondHLCPNm\lnP{:=}\prcondHLCPNm#3\\#4&\prcondHLCPNm\lnP{:=}\prcondHLCPNm#5, \end{split}\end{equation} }
		\newcommand{\envHLineCFDefNmPd}[5][*]{ \begin{equation}\begin{split} \ifx*#1 \text{[ERROR;need label name]} \else \label{#1} \fi #2&\prcondHLCPNm\lnP{:=}\prcondHLCPNm#3,\\#4&\prcondHLCPNm\lnP{:=}\prcondHLCPNm#5. \end{split}\end{equation} }
		\newcommand{\envHLineCSCmNm}[7][*]{ \begin{equation}\begin{split} \ifx*#1 \text{[ERROR;need label name]} \else \label{#1} \fi #2&\prcondHLCPNm\lnP{=}\prcondHLCPNm#3,\\#4&\prcondHLCPNm\lnP{=}\prcondHLCPNm#5,\\#6&\prcondHLCPNm\lnP{=}\prcondHLCPNm#7 \end{split}\end{equation} }
		\newcommand{\envHLineCSNm}[7][*]{ \begin{equation}\begin{split} \ifx*#1 \text{[ERROR;need label name]} \else \label{#1} \fi #2&\prcondHLCPNm\lnP{=}\prcondHLCPNm#3\\#4&\prcondHLCPNm\lnP{=}\prcondHLCPNm#5\\#6&\prcondHLCPNm\lnP{=}\prcondHLCPNm#7 \end{split}\end{equation} }
		\newcommand{\envHLineCSNmPd}[7][*]{ \begin{equation}\begin{split} \ifx*#1 \text{[ERROR;need label name]} \else \label{#1} \fi #2&\prcondHLCPNm\lnP{=}\prcondHLCPNm#3,\\#4&\prcondHLCPNm\lnP{=}\prcondHLCPNm#5,\\#6&\prcondHLCPNm\lnP{=}\prcondHLCPNm#7. \end{split}\end{equation} }
		\newcommand{\envHLineCSCmDefNm}[7][*]{ \begin{equation}\begin{split} \ifx*#1 \text{[ERROR;need label name]} \else \label{#1} \fi #2&\prcondHLCPNm\lnP{:=}\prcondHLCPNm#3,\\#4&\prcondHLCPNm\lnP{:=}\prcondHLCPNm#5,\\#6&\prcondHLCPNm\lnP{:=}\prcondHLCPNm#7 \end{split}\end{equation} }
		\newcommand{\envHLineCSDefNm}[7][*]{ \begin{equation}\begin{split} \ifx*#1 \text{[ERROR;need label name]} \else \label{#1} \fi #2&\prcondHLCPNm\lnP{:=}\prcondHLCPNm#3\\#4&\prcondHLCPNm\lnP{:=}\prcondHLCPNm#5\\#6&\prcondHLCPNm\lnP{:=}\prcondHLCPNm#7 \end{split}\end{equation} }
		\newcommand{\envHLineCSDefNmPd}[7][*]{ \begin{equation}\begin{split} \ifx*#1 \text{[ERROR;need label name]} \else \label{#1} \fi #2&\prcondHLCPNm\lnP{:=}\prcondHLCPNm#3,\\#4&\prcondHLCPNm\lnP{:=}\prcondHLCPNm#5,\\#6&\prcondHLCPNm\lnP{:=}\prcondHLCPNm#7. \end{split}\end{equation} }
		\newcommand{\envHLineCECmNm}[9][*]{ \begin{equation}\begin{split} \ifx*#1 \text{[ERROR;need label name]} \else \label{#1} \fi #2&\prcondHLCPNm\lnP{=}\prcondHLCPNm#3,\\#4&\prcondHLCPNm\lnP{=}\prcondHLCPNm#5,\\#6&\prcondHLCPNm\lnP{=}\prcondHLCPNm#7,\\#8&\prcondHLCPNm\lnP{=}\prcondHLCPNm#9,  \end{split}\end{equation} }
		\newcommand{\envHLineCENm}[9][*]{ \begin{equation}\begin{split} \ifx*#1 \text{[ERROR;need label name]} \else \label{#1} \fi #2&\prcondHLCPNm\lnP{=}\prcondHLCPNm#3\\#4&\prcondHLCPNm\lnP{=}\prcondHLCPNm#5\\#6&\prcondHLCPNm\lnP{=}\prcondHLCPNm#7\\#8&\prcondHLCPNm\lnP{=}\prcondHLCPNm#9  \end{split}\end{equation} }
		\newcommand{\envHLineCENmPd}[9][*]{ \begin{equation}\begin{split} \ifx*#1 \text{[ERROR;need label name]} \else \label{#1} \fi #2&\prcondHLCPNm\lnP{=}\prcondHLCPNm#3,\\#4&\prcondHLCPNm\lnP{=}\prcondHLCPNm#5,\\#6&\prcondHLCPNm\lnP{=}\prcondHLCPNm#7,\\#8&\prcondHLCPNm\lnP{=}\prcondHLCPNm#9.  \end{split}\end{equation} }
		\newcommand{\envHLineCECmDefNm}[9][*]{ \begin{equation}\begin{split} \ifx*#1 \text{[ERROR;need label name]} \else \label{#1} \fi #2&\prcondHLCPNm\lnP{:=}\prcondHLCPNm#3,\\#4&\prcondHLCPNm\lnP{:=}\prcondHLCPNm#5,\\#6&\prcondHLCPNm\lnP{:=}\prcondHLCPNm#7,\\#8&\prcondHLCPNm\lnP{:=}\prcondHLCPNm#9,  \end{split}\end{equation} }
		\newcommand{\envHLineCEDefNm}[9][*]{ \begin{equation}\begin{split} \ifx*#1 \text{[ERROR;need label name]} \else \label{#1} \fi #2&\prcondHLCPNm\lnP{:=}\prcondHLCPNm#3\\#4&\prcondHLCPNm\lnP{:=}\prcondHLCPNm#5\\#6&\prcondHLCPNm\lnP{:=}\prcondHLCPNm#7\\#8&\prcondHLCPNm\lnP{:=}\prcondHLCPNm#9  \end{split}\end{equation} }
		\newcommand{\envHLineCEDefNmPd}[9][*]{ \begin{equation}\begin{split} \ifx*#1 \text{[ERROR;need label name]} \else \label{#1} \fi #2&\prcondHLCPNm\lnP{:=}\prcondHLCPNm#3,\\#4&\prcondHLCPNm\lnP{:=}\prcondHLCPNm#5,\\#6&\prcondHLCPNm\lnP{:=}\prcondHLCPNm#7,\\#8&\prcondHLCPNm\lnP{:=}\prcondHLCPNm#9.  \end{split}\end{equation} }
	\newcommand{\envMLineTPd}[3][*]{ \ifx*#1 \begin{multline*} #2\lnP{=}#3.\end{multline*} \else \begin{multline} \label{#1} #2\lnP{=}#3.\end{multline} \fi }
	\newcommand{\envMLineTCm}[3][*]{ \ifx*#1 \begin{multline*} #2\lnP{=}#3,\end{multline*} \else \begin{multline} \label{#1} #2\lnP{=}#3,\end{multline} \fi }
	\newcommand{\envMLineTDefPd}[3][*]{ \ifx*#1 \begin{multline*} #2\lnP{:=}#3.\end{multline*} \else \begin{multline} \label{#1} #2\lnP{:=}#3.\end{multline} \fi }
	\newcommand{\envMLineTCmDef}[3][*]{ \ifx*#1 \begin{multline*} #2\lnP{:=}#3,\end{multline*} \else \begin{multline} \label{#1} #2\lnP{:=}#3,\end{multline} \fi }

	\newcommand{\envCaseTCm}[3][?]{\begin{cases} \glcondDis{#1}#2,\lcparaCase\\\glcondDis{#1}#3,\end{cases}}
	\newcommand{\envCaseTPd}[3][?]{\begin{cases} \glcondDis{#1}#2,\lcparaCase\\\glcondDis{#1}#3.\end{cases}}
	\newcommand{\envCaseThCm}[4][?]{\begin{cases} \glcondDis{#1}#2,\lcparaCase\\\glcondDis{#1}#3,\lcparaCase\\\glcondDis{#1}#4,\end{cases}}

	\newcommand{\ldgP}[2][]{\Fc[#1]{\deg\,}{#2}}
		\newcommand{\dgP}{\ldgP}		

	\newcommand{\lfsPT}[3][n]{\mathcal{P}_{#2}^{#3}}
	\newcommand{\lfsQO}[2][n]{\ifx m#1\mathcal{\mLt{Q}}^{#2}\else\mathcal{Q}^{#2}\fi}
	\newcommand{\lfsQT}[3][n]{\ifx m#1\mathcal{\mLt{Q}}_{#2}^{#3}\else\mathcal{Q}_{#2}^{#3}\fi}

				\newcommand{\fsP}{\lfsPT}
				
				\newcommand{\fsQ}{\lfsQT}

	\newcommand{\lfmP}[2][n]{\mathbf{P}_{#2}}
	\newcommand{\lfmQ}[2][n]{\ifx m#1\mathbf{\mLt{Q}}_{#2}\else\mathbf{Q}_{#2}\fi}
	
	\newcommand{\lfmX}{{\bf X}}

		\newcommand{\fmP}{\lfmP}	
		\newcommand{\fmQ}{\lfmQ}		
		
		\newcommand{\fmX}{\lfmX}	

	\newcommand{\lvlT}{\mathrm{True}}	
	\newcommand{\lvlF}{\mathrm{False}}	
		\newcommand{\vlT}{\lvlT}
		\newcommand{\vlF}{\lvlF}

	\newcommand{\lasr}[2][]{\ifx a#2{\rm(a)}\else\ifx b#2{\rm(b)}\else??\fi\fi }	
		
	\newcommand{\ltecO}[2][?]{({\bf T#2})}
		\newcommand{\tec}{\ltecO}
		
	\newcommand{\lprp}[2][]{\ifx1#2{\rm(A)}\else\ifx2#2{\rm(Z?)}\else\ifx3#2{\rm(B)}\else\ifx4#2{\rm(C)}\else??\fi\fi\fi\fi }	
	\newcommand{\lprpp}[2][]{\ifx1#2{\rm(A)$\mbox{}_{\mLt{Q}}$}\else\ifx2#2{\rm(Z?)$\mbox{}_{\mLt{Q}}$}\else\ifx3#2{\rm(B)$\mbox{}_{\mLt{Q}}$}\else\ifx4#2{\rm(C)$\mbox{}_{\mLt{Q}}$}\else??\fi\fi\fi\fi }
	\newcommand{\lprpS}[2][]{\ifx1#2{\rm(A)$\mbox{}_{Q}$}\else\ifx2#2{\rm(Z?)$\mbox{}_{Q}$}\else\ifx3#2{\rm(B)$\mbox{}_{Q}$}\else\ifx4#2{\rm(C)$\mbox{}_{Q}$}\else??\fi\fi\fi\fi }
		\newcommand{\prp}{\lprp}
		\newcommand{\prpp}{\lprpp}
		\newcommand{\prpS}{\lprpS}

	\newcommand{\ltplE}{\varnothing}		
		\newcommand{\tplE}{\ltplE}

	\newcommand{\llaPO}[2][n]{\Fc[#1]{\pi}{#2}}		
	\newcommand{\llaPT}[3][n]{\idFc[#1]{\pi}{#2}{#3}}		

	\newcommand{\lsbF}[2][?]{\idFc{\llaPO[]{}}{\mathrm{min}}{#2}}		
	\newcommand{\lsbE}[2][?]{\idFc{\llaPO[]{}}{\mathrm{max}}{#2}}		
		\newcommand{\sbF}{\lsbF}
		\newcommand{\sbE}{\lsbE}

	\newcommand{\lclG}[3][n]{\idFc[#1]{g}{#2}{#3}}	
		\newcommand{\clG}{\lclG}

	\newcommand{\lplP}[3][n]{\idFc[#1]{P}{#2}{#3}}		
		\newcommand{\plP}{\lplP}
	\newcommand{\lplF}[3][n]{\idFc[#1]{F}{#2}{#3}}		
		\newcommand{\plF}{\lplF}

	\newcommand{\lclES}[1][?]{\mathsf{ESat}}	

	\newcommand{\lclSN}{\mathsf{Sat}}	
	\newcommand{\lclSO}[1]{\lclSN_{#1}}	
	\newcommand{\lclST}[2]{\lclSN_{#1}^{(#2)}}	
				\newcommand{\clSO}{\lclSO}		\newcommand{\clST}{\lclST}

	\newcommand{\lrgBP}{\mathbb{BP}}	
		\newcommand{\rgBP}{\lrgBP}

	\newcommand{\lidI}[1][?]{\ifx #1? \mathfrak{I} \else\mathfrak{I} _#1\fi}	
		\newcommand{\idI}{\lidI}

	\newcommand{\llxor}{\oplus}
		\newcommand{\lxor}{\llxor}

	\newcommand{\lclM}[2][n]{\mLt{c}}

	\newcommand{\lvrsMN}[1][n]{{x}_{\mathrm{min}}}	
	\newcommand{\lvrJ}[1][?]{x_{wild}}		

	\newcommand{\lstL}{\mathcal{L}}		

		\newcommand{\stL}{\lstL}

	\newcommand{\lbFc}[3][n]{\Fc[#1]{\eLt{#2}}{#3}}	
		\newcommand{\bFc}{\lbFc}

	\newcommand{\lspB}[1][?]{\mathcal{BF}}	
	\newcommand{\lspCLO}[1]{\mathbb{CL}_{#1}}
	\newcommand{\lspCLT}[2]{\mathbb{CL}_{#1}^{(#2)}}
		
		\newcommand{\spB}{\lspB}
				\newcommand{\spCLO}{\lspCLO}		\newcommand{\spCLT}{\lspCLT}

	\newcommand{\lsyN}{\mathrm{S}}	
	\newcommand{\lsy}[2][n]{\Fc[#1]{\lsyN}{#2}}	
		\newcommand{\syN}{\lsyN}
		\newcommand{\sy}{\lsy}
	\newcommand{\lsyNN}{\mathrm{T}}	
		\newcommand{\syNN}{\lsyNN}
		
	\newcommand{\leqV}{\sim}	
	\newcommand{\leqS}{\approx}	
	\newcommand{\leqE}[1][?]{\,\ifx ?#1\leqV\else{\leqV}_{#1}\fi\,}		
		
		\newcommand{\eqS}{\leqS}

	\newcommand{\ldlC}[2][n]{#2'}		
		\newcommand{\dlC}{\ldlC}
		
	\newcommand{\lmpF}[2][n]{\Fc[#1]{\varphi}{#2}}		
	\newcommand{\lmpFi}[2][n]{\pwFc[#1]{\varphi}{-1}{#2}}		
	\newcommand{\lmpN}[2][n]{\Fc[#1]{\mathcal{N}}{#2}}		
	\newcommand{\lmpE}[2][n]{\Fc[#1]{\mathcal{\psi}}{#2}}		

		\newcommand{\mpF}{\lmpF}	\newcommand{\mpFi}{\lmpFi}
		\newcommand{\mpN}{\lmpN}
		
		\newcommand{\mpE}{\lmpE}

	\newcommand{\lpaBW}[1][o]{B}			
	\newcommand{\lpaTW}[1][m]{t}		
		\newcommand{\paBW}{\lpaBW}
		\newcommand{\paTW}{\lpaTW}

\allowdisplaybreaks[4]
	\title{\mainTitle}
	\author{\authorName
			\thanks{\organizationNameFst, \placeAddressFst}
		}
	\date{}

\begin{document}
\maketitle
\renewcommand{\thefootnote}{\fnsymbol{footnote}}
\footnote[0]{e-mail : \emailAddressFst}
\footnote[0]{MSC-class: \MSCname}
\renewcommand{\thefootnote}{\arabic{footnote}}\setcounter{footnote}{0}
\vPack[30]

\begin{abstract}
It is known a method for converting  
	a system of Boolean polynomial equations to a single Boolean polynomial equation with less variables.
In this paper,
	we 
	show a formula for systems of Boolean polynomial equations
	which is based on the method.
The formula 		
	has a structure of binary tree,
	and
	conforms to De Morgan's duality.
Using the formula,
	we prove a computational complexity result 
	with
	a parameter
	for solving systems.
The parameter 
	is	
	the bandwidth in matrix and graph theories:
	to be precise,	
	the definition follows convention in matrix
	and
	the value depends on the order of variables.	
We also apply the result 
	to
	the NP-complete problems, 
	SAT and graph list-coloring,
	to show that
	these problems are fixed parameter tractable by bandwidth.
	
\end{abstract}

\section{Introduction} \label{sectOne}
The finite field $\setF[2]=\SetO{0,1}$ with two elements,
	which is also called the Galois field $\mathbb{GF}(2)$ in his honor, 
	plays fundamental roles in mathematics and computer science.		
It is the smallest finite field 
	and
	its algebraic rules 	
	are 
	determined by a few equations involving the addition ``$+$'' and multiplication ``$\,\cdot\,$''.
One of the outstanding facts of $\setF[2]$
	is a structural relation 
	to the two-element Boolean algebra $\setB=\SetO{\vlF,\vlT}$
	under the identifications of $\vlF=0$ and $\vlT=1$.	
That is,
	for any pair $(\alp, \beta)$ of elements, 	
	\envPLineCm[1_PL_RelOpeBolF2]
	{
		\alp \land \beta	 	\lnP{=}	\alp \cdot \beta
	,\qquad
		\alp \lor \beta		\lnP{=}	(\alp+1) \cdot (\beta+1) + 1
	,\qquad
		\alp \lxor \beta 		\lnP{=}	\alp + \beta
	}
	where 
	$\land$, $\lor$, and $\lxor$ 
	stand for 
	the binary operations of conjunction, disjunction, and exclusive disjunction in $\mathbb{B}$, 
	respectively.
The unary operation $\lnot$ of negation is expressed as $\lnot \alp = \alp + 1$.	

A Boolean polynomial,
	which is also called a Boolean expression in algebraic normal form \cite{CG99},	
	Reed-Muller expansion \cite{Muller54, Reed54},
	and
	Zhegalkin polynomial \cite{Gegalkine27},
	naturally arises 		
	when we transform  a Boolean expression  
	to a polynomial 	
	using \refEq{1_PL_RelOpeBolF2}.	
The polynomial 
	is 
	a congruence class of the polynomial ring $\setF[2][x_1, \ldots, x_n]$ in $n$ variables, 
	and
	identified with a Boolean function from $\setF[2]^n$ to $\setF[2]$.
(Details will be introduced in \refSect{sectTwo}.)
The Boolean polynomials and the ring consisting of them
	are important subjects in various areas:   		
	e.g., 
	algebraic geometry \cite{BIKPP96, CLO15, Kunz85}, 
	Boolean ideal and variety \cite{Lundqvist15,RM14},
	circuit theory \cite{Williams14},
	cording theory \cite{GOS10,Meneghetti18}, 
	cryptography \cite{CG99,JV18}, 
	and
	Gr\"obner basis \cite{BD09, CEI96, SISNKS11}.	
Although the contexts differ depending on the areas,
	solving a system of Boolean polynomial equations is a common problem.
	
Recently,
	\etalTx{Lokshtanov} \cite{LPTWY17} used several techniques developed from circuit complexity
	to construct algorithms 
	for the problem,		
	which  
	beat brute force search 
	without relying on any heuristic conjectures.
(They actually studied for not only $\setF[2]$ but also any finite fields.)
In this paper, 	
	we focus on the two basic techniques in \cite{LPTWY17}:
	\tec{1} transform a system of Boolean polynomial equations to a single Boolean polynomial equation;
	and
	\tec{2} transform a single Boolean polynomial equation to one with less variables.
It may be worth noting that 
	\tec{1} is a classical fact in algebraic geometry; 
	for example,
	see Exercise 3 in \cite[Chapter I, Section 1]{Kunz85}.
Combining \tec{1} and \tec{2}, 	
	we can convert a system of Boolean polynomials to a single Boolean polynomial equation with less variables.

The aims of this paper are
	to show a formula 
	which is based on the converting method, 	
	and
	to give applications to computational complexity.	
Considering the method from a sequential viewpoint,			
	we construct 
	a formula for systems of Boolean polynomial equations  (\refThm{3_Thm1}).	
Using the formula,
	we prove 
	a parameterized complexity result for solving systems (\refThm{3_Thm2}).		
Then we apply the complexity result 
	to NP-complete problems: SAT and graph list-coloring (\refCor[s]{3_Cor1} and \ref{3_Cor2}).
	
Sketches of our results are following.
The details will be stated in \refSect{sectThree}.		

The formula of \refThm{3_Thm1}
	possesses	
	both operations of conjunction and disjunction recursively,
	with a structure of binary tree.
By the recursiveness and structure,
	the formula
	conforms to De Morgan's duality.
The distributivity of the operations 
	plays a fundamental role in the proof.

\refThm{3_Thm2} follows from 
	the fact that
	it is possible 
	to reduce leaf nodes on the binary tree in some cases.
Let
	$\sbL[]{}$ be the big O notation,	
	and
	let $\sbLp[]{}$ denote the notation which omits polynomial factors in $\sbL[]{}$.	
\refThm{3_Thm2} implies that	
	the satisfiability of a system is decidable in time $\sbLp{2^{\paBW}}$,
	where
	$\paBW$ is 	
	the bandwidth in matrix and graph theories.	
The definition of bandwidth in this paper will adopt convention in matrix;
	that is,
	the value is not minimum
	but 
	depends on the order of variables.		
The annihilator and identity laws for disjunction are crucial in the proof.
We note that		
	polynomial factors are not omitted in the actual statement.		

\refCor[s]{3_Cor1} and \ref{3_Cor2} are consequences of 
	the fact that 
	systems of Boolean polynomial equations can express  
	the NP-complete problems,
	CNF-SAT, BMQ-SAT and graph list-coloring.
The CNF-SAT problem,
	the Boolean satisfiability problem in conjunctive normal form,	
	is
	the first NP-complete problem \cite{Cook71, Levin73}.
CNF-SAT 
	has many applications in the real world \cite{BHMW09}. 
The BMQ-SAT problem,
	for which 
	the algorithm beating brute force search was presented in \cite{LPTWY17},	
	is 
	the satisfiability problem of a Boolean multivariate quadratic system 
	(or a system of Boolean polynomial equations of degree $2$).	
BMQ-SAT is significant in cryptography 	
	to generate secure ciphers \cite{BCJ07, BFSS13}.
We mean by SAT either one of both.
The graph list-coloring problem
	is 
	a generalization of the original coloring problem:
	in addition to the proper condition such that no two adjacent vertices receive the same color,
	a list of allowed colors is imposed for each vertex.
Graph coloring is a central problem as SAT is, in theoretical, practical and historical aspects \cite{GJPS17}.

A problem of input size $n$ with a parameter $k$ 	
	is called 
	fixed parameter tractable (or FPT for short)
	if it can be solved in time $\Fc{f}{k} n^{\sbL{1}}$, 
	where $f$ is a function only depending on $k$.
Parameterized complexity theory
	is 
	a two dimensional analog of the classical framework of P versus NP,		
	and
	class FPT corresponds to class P
	(see \cite{FFLRSST11} and references therein for details).
Our complexity results show that
	SAT and list-coloring are FPT by bandwidth $\paBW$,
	where,
	in the latter problem,
	the total number $l$ of allowed colors 
	is considered to be 
	constant and independent to the size $n$. 	
	
It appears that
	our complexity result of list-coloring problem
	is especially interesting,
	because  
	the problem is known to be
	W[1]-hard 
	for both parameters of treewidth and vertex cover \cite{FFLRSST11,FLMRS08,FGK11},
	where
	W[1] is the class corresponding to NP.
That is,
	the bandwidth 
	is 	
	a different type parameter in list-coloring.
The function $f$ is roughly $(2l)^{\paBW}$.
Our complexity result of CNF-SAT is already known (in a sense),
	because
	CNF-SAT is FPT by treewidth (of incidence graph) \cite{Szeider03}
	and
	treewidth is more general than bandwidth.
However
	our result has an advantage:
	$f$ is expressed as $2^{\paBW}$ 	
	and
	it is concrete;		
	in contrast, 
	$f$ in \cite{Szeider03} is abstract (see Theorem 4 and Corollary 1 in the paper).
Our complexity result of BMQ-SAT	
	seems to be new.	

The problem of finding bandwidth is NP-hard
	and
	its decisional version is NP-complete.
However 
	there are 
	many heuristic algorithms including the Cuthill-McKee algorithm,
	and
	polynomial-time algorithms for spacial classes of graphs.	
(See \cite{CCDE82,DPS02} and references therein for details.)
Randomized approximate algorithms for general graphs,
	which run 
	in polynomial or nearly linear time
	and
	have polylogarithmic factors of optimal,
	are also known \cite{BKRV98,Feige00}.
Thanks to those algorithms,
	our complexity results 
	are
	practical
	if the bandwidth $\paBW$ is small.

The paper is organized as follows.
In \refSec{sectTwo},			
	we quickly review the Boolean polynomials and their basic properties.	
Rigid statements of our results are given in \refSec{sectThree}.	
We prove  
	\refThm{3_Thm1} in \refSec{sectFour}, 		
	and 
	\refThm{3_Thm2} in \refSec{sectFive}.
\refSec{sectSix}
	is devoted 
	to the proofs of \refCor[s]{3_Cor1} and \ref{3_Cor2}.

\section{Review of the Boolean polynomials} \label{sectTwo}	
The finite field $\setF[2]$ is commutative,
	and
	its algebraic rules are determined 	
	by the equations involving the addition and multiplication: 
	\envPLinePd
	{
		0+0=1+1=0\cdot0=0\cdot1=0
	,\qquad
		0+1=1\cdot1=1
	}
The subtraction and division are unnecessary,
	because
	the subtraction is identical to the addition
	and
	no invertible elements except $1$ exist.

The Boolean polynomial ring $\rgBP_n=\rgBP[x_1, \ldots, x_n]$ 
	is defined by the quotient ring
	\envHLineCmDef[3_PL_DefBooPol]
	{
		\rgBP_n	
	}
	{
		\setF[2][x_1, \ldots, x_n] / \idI[n]	
	}
	where 	
	\envHLineDefPd
	{
		\idI[n]	
	}
	{
		\SetT{ p_1(x_1^2+x_1) + \cdots + p_n(x_n^2+x_n) }{ p_i \in \setF[2][x_1, \ldots, x_n] }
	}  	
A Boolean polynomial $p=\Fc{p}{x_1, \ldots, x_n}$ is a congruence class in $\rgBP_n$.	
In the ring,
	the variables are idempotent (i.e., $x_j^2=x_j$),
	and
	the number of monomials is $2^n$.
Because 
	the monomials are independent,
	$p$ is uniquely expressed as		
	\envPLinePd[3_PL_ExpBolPol]
	{
		p			
	\lnP{=}	
		\Sm{e_1, \ldots, e_n \in \SetO{0,1}} a_{e_1, \ldots, e_n}x_{1}^{e_1} \cdots x_{n}^{e_n}
	\qquad
		(a_{e_1, \ldots, e_n} \in \setF[2])
	}	
We thus have  
	\envMCm
	{
		\nmSet{\rgBP_n} 
	}
	{ 
		2^{2^n}
	}	
	where
	we mean by $\nmSet{A}$ 
	the number of elements of a set $A$.	

Let $\spB_n$ be the ring of Boolean functions of $n$ variables,
	or the ring of $\setF[2]$-valued functions with the domain $\setF[2]^n$.
For a Boolean polynomial $p=\Fc{p}{x_1, \ldots, x_n}$,
	we denote by $\bFc[]{p}{}$ the polynomial function of $p$,
	which 	
	is defined by 	
	\envPLine
	{\begin{array}{cccc}
		\bFc[]{p}{}:	&\setF[2]^n 			&\quad \to \quad& \setF[2] \\
				&\roTx{90}{$\in$} 		&  		& \roTx{90}{$\in$} \\
				&(\alp_1, \ldots, \alp_n) 	& \mapsto & \Fc{p}{\alp_1, \ldots, \alp_n}.
	\end{array}}
This induces a well-defined homomorphism from $\rgBP_n$ to $\spB_n$,
	since 
	$0^2+0=1^2+1=0$		
	and
	$\bFc[]{p}{}$ is the zero function if $p$ is in $\idI[n]$.		
The homomorphism		
	is isomorphic,\footnote{
For the injectivity,
	we may show that $\bFc[]{p}{}\neq 0$ for a non-zero Boolean polynomial $p$,
	which follows from the unique expression in \refEq{3_PL_ExpBolPol}. 	
For the surjectivity, 
	we may show that the numbers of elements in both rings are equal,
	or $\nmSet{\spB_n}=2^{2^n}$,		
	which 
	follows from
	$\nmSet{\setF[2]^n} = 2^n$ and $\nmSet{\setF[2]}=2$.
}
	and 
	$\rgBP_n$ can be identified with $\spB_n$: 	
	\envPLine[2_PL_IdtPolFnc]
	{\begin{array}{ccc}
		\rgBP_n 			&\quad \simeq \quad& \spB_n \\
		\roTx{90}{$\in$} 		&  		& \roTx{90}{$\in$} \\
		p	& \leftrightarrow & \,\bFc[]{p}{}\,.
	\end{array}}

We see from \refEq{2_PL_IdtPolFnc}
	that
	$\rgBP_n$
	has
	the same calculation rules as the codomain of $\spB_n$, or each of $\setF[2]$ and $\setB$.
Thus
	the identities in \refEq{1_PL_RelOpeBolF2} hold on $\rgBP_n$.
In addition,
	we have
	modular arithmetic properties 
	\envPLineCm[2_PL_EqBasBooPolF]
	{
		p + p		\lnP{=}	0		
	,\qquad	
		p^2 		\lnP{=}	p
	}
	and
	annihilator and identity laws
	\envPLinePd[2_PL_EqBasBooPolB1]
	{
		0 \land p			\lnP{=}	0
	,\qquad
		1 \lor p			\lnP{=}	1
	,\qquad	
		1 \land p			\lnP{=}	0 \lor p			\lnP{=}	p
	}
Generalizing the second equation in \refEq{1_PL_RelOpeBolF2} to $m$ elements,
	we also have
	\envHLineCm[2_PL_RelAndOr]
	{
		p_1 \lor \cdots \lor p_m
	}
	{
		(p_1+1) \cdots (p_m+1) + 1
	}
	where
	$p_1, \ldots, p_m$ are Boolean polynomials.
(For \refEq{2_PL_RelAndOr}, 
	see, e.g., 
	\cite[Section 3]{Williams14}.)

\section{Statement of results} \label{sectThree}	
We begin with preparing notations and terminologies.	

For a pair $(\syN, \syNN)$ of systems,
	we say
	that
	$\syN \eqS \syNN$ 
	if either both systems are satisfiable or both are not.
It is easily seen that
	$\eqS$ is an equivalence relation.
We call a system including two or more equations a multiple system;
	in contrast,
	we call a system including only one equation a single system.
A system means either one of both.	
Let $i$ be a positive integer at most $n+1$.
For a positve integer $j$ less than $n+1$,
	we denote by $\rgBP_{i,j}$ 
	the subring $\rgBP[x_{i}, \ldots, x_{j}]$ of $\rgBP_{n}$, 
	where
	$\rgBP_{i,j} = \rgBP[\phi] = \setF[2]$
	if
	$i > j$. 
It holds that
	\envHLinePt[3_PL_EqBooPolSub]{\in}
	{
		p |_{x_{i}=\alp} 
	}
	{ 
		\rgBP_{i+1,j} 		
	} 
	when $p \in \rgBP_{i,j}$ and $\alp \in \setF[2]$.	

We will describe the two basic techniques in \cite{LPTWY17},
	which enable us to 		
	convert 
	a multiple system to a single system with less variables.
Let $\syN=\sy{p_1, \ldots, p_m}$ be a system of Boolean polynomials $p_1, \ldots, p_m$. 
We define a Boolean polynomial 
	by 
	\envHLineDefPd[3_PL_DefPofSys]
	{	
		\plP[]{\syN}{}
	}
	{
		p_1 \lor \cdots \lor p_m
	\lnP{=}
		(p_1+1) \cdots (p_m+1) + 1
	}	
The annihilator and identity laws for $\lor$ 
	imply that
	$\plP[]{\syN}{} = 0$
	if and only if 
	$p_1 = \cdots = p_m = 0$.
Hence
	$\syN$ 
	is equivalent 
	to the single system consisting of $\plP[]{\syN}{}$ under $\eqS$,
	and
	we can apply $\plP[]{\syN}{}$ to solve $\syN$.
This is one of the techniques, 
	which enable us to transform a multiple system to a single system.
We put
	$\plP[]{1}{} = \plP[]{\syN}{}$.
For an integer $j$ from $2$ to $n+1$,	
	we recursively define 
	a Boolean polynomial in $\rgBP_{j,n}$ by	
	\envHLineDefPd[3_PL_DefPolP]
	{
		\plP[]{j}{} 	
	}
	{
		\Pd{\alp \in \setF[2]}	\bkR{ \plP[]{j-1}{} |_{x_{j-1}=\alp} } 
	} 		
The number of variables in $\plP[]{j}{}$ is at most $n+1-j$,
	and
	it decreases as $j$ increases.
Let $\syN_j$ be a single system consisting of $\plP[]{j}{} $ for each $j$.
Obviously,
	$\syN_{j-1} \eqS \syN_{j}$,
	and
	we can apply $\plP[]{j}{}$ to solve $\plP[]{1}{} = \plP[]{\syN}{}$.	
This is another technique to reduce variables.		
Combining these techniques,
	we can use 
	$\plP[]{j}{} \in \rgBP_{j,n}$ to solve $\syN$.			
	
More notations will be required
	to state \refThm{3_Thm1}.
We will first 
	introduce the definition of the CNF-SAT problem,
	next 
	define notations on systems which involve CNF-SAT,
	and
	then mention the others.

Let $x$ be a variable.
To distinguish $x$ and $\lnot x = x+1$,	
	we call the former a positive literal
	and
	the latter a negative literal.
A literal means either one of both.
The CNF-SAT problem
	is the problem of 
	deciding if there exists an assignment of variables 
	which satisfies a conjunction of clauses,	
	where 
	a clause means a disjunction of literals.
For instance,
	a CNF-SAT problem
	is
	solving	
	\envHLineCm[3_PL_EqBolEx]
	{
		x_1	\land		(\lnot x_1 \lor	x_2)		
	}	
	{
		\vlT
	}
	which is satisfiable 
	because 
	\envM
	{
		(x_1, x_2)
	}
	{
		(\vlT, \vlT)
	}
	is a solution.
It is easily seen 
	from 		
	the annihilator and identity laws for $\land$	
	that
	\refEq{3_PL_EqBolEx} 
	is equivalent to the system of Boolean equations,
	\envPLine[3_PL_EqBolExSys]
	{
		\left\{\begin{array}{ccl}
			x_1 				& = & 1, \\
			\lnot x_1 \lor x_2 	& = & 1.
		\end{array}\right.
	}	
	
We define a subspace in $\rgBP_n$ by 
	\envHLineCmDef[3_PL_DefCls]
	{
		\spCLO{n}
	}
	{
		\setF[2] \cup \SetT{  l_1 \cdots l_k }{ k\geq1, \text{$l_i$ are literals} }		
	}
	and
	its extension by
	\envHLineDefPd[3_PL_DefClsR]
	{
		\spCLT{n}{j}
	}
	{
		\usbTx{j}{\spCLO{n} + \cdots + \spCLO{n}}
	\lnP{=}
		\SetT{ c_1 + \cdots + c_j }{  c_1, \ldots, c_j \in \spCLO{n} }
	}
We have 
	$\spCLT{n}{2^n} = \rgBP_n$
	since
	$\spCLO{n}$ includes all monomials in $\rgBP_n$.
Let $c=l_1 \cdots l_k$ be a non-constant Boolean polynomial in $\spCLO{n}$.
When $l_i = l_j$,	
	$l_i l_j = l_i = l_j$ by \refEq{2_PL_EqBasBooPolF}
	and we can remove either $l_i$ or $l_j$ from $c$.
When $l_i = l_j+1$,		
	$l_i l_j = 0$ by \refEq{2_PL_EqBasBooPolF}
	and $c$ is the zero polynomial,
	which contradicts the non-constant.
Therefore,		
	in this paper,			
	we will assume that		
	the literals $l_1, \ldots, l_k$ appearing in a polynomial of $\spCLO{n}$ satisfy
	\envPLinePd[2_PL_AsmCls]			
	{	
		l_i \lnP{\notin} \Set{l_j, l_j+1} 	
		\quad\text{for}\quad 
		i	\lnP{\neq} j
	}	
For each literal $l_i$,
	let $y_i$ and $\alp_i$ 
	denote a variable in $\Set{x_1, \ldots, x_n}$
	and 
	a value in $\setF[2]$,
	respectively,
	such that
	$l_i = y_i + \alp_i$.	
We have
	the following correspondence 	
	between 		
	equations of a polynomial and a clause:\footnote{
The following equivalences hold by \refEq{1_PL_RelOpeBolF2} and De Morgan's duality:
	\envMFPdPt{\Leftrightarrow}
	{
		(y_1 + \alp_1)  \cdots  (y_k + \alp_k)  = 0
	}
	{
		(y_1 + \alp_1)  \land \cdots \land  (y_k + \alp_k)  = \vlF
	}
	{
		(y_1 + \alp_1 + 1)  \lor \cdots \lor  (y_k + \alp_k + 1)  = \vlT
	}
	{
		\bullet_1 y_1 \lor \cdots \lor \bullet_k y_k			=	\vlT
	}
}
	\envPLineCm[3_PL_CorPolCls]
	{
		(y_1 + \alp_1)  \cdots  (y_k + \alp_k)  		\lnP{=}	0
	\qquad\Leftrightarrow\qquad
		\bullet_1 y_1 \lor \cdots \lor \bullet_k y_k	\lnP{=}	\vlT
	}
	where 
	$\bullet_i$ stands for the negation `$\lnot$' if $\alp_i=0$
	and 
	the empty letter if $\alp_i=1$. 
For instans,
	$y_1+1=0$ and $y_1 (y_2+1) = 0$ correspond to		
	$y_1 = \vlT$ and $\lnot y_1 \lor y_2 = \vlT$,	
	respectively.
Therefore we call an element of $\spCLO{n}$ a clause polynomial, 
	or 
	simply a clause.	
Because of \refEq{3_PL_CorPolCls}
	and 
	the correspondence 
	between \refEq{3_PL_EqBolEx} and \refEq{3_PL_EqBolExSys},		
	the set of CNF-SAT problems in $n$ variables
	is equivalent to 		
	\envHLineDefPd[3_PL_DefSysSAT]
	{
		\clSO{n}	
	}	
	{
		\SetT{ \sy{c_1, \ldots, c_m} }{ c_1, \ldots, c_m \in \spCLO{n}}	
	}
As an extension of \refEq{3_PL_DefSysSAT},	
	we define		
	\envHLineCmDef[3_PL_DefSysSATSum]
	{
		\clST{n}{l}	
	}	
	{
		\SetT{ \sy{c_1, \ldots, c_m} }{ c_1, \ldots, c_m \in \spCLT{n}{l} }	
	}
	where 
	$l$ is a positive integer.
Since
	$\spCLT{n}{2^n} = \rgBP_n$,	
	$\clST{n}{2^n}$ covers all systems of Boolean polynomial equations.

For a system $\syN=\sy{p_1, \ldots, p_m}$,
	we call $k=\tpMax{ i }{ \dgP{p_i} }$ the degree of $\syN$;
	the system is usually called a $k$-CNF-SAT problem 	
	if 
	$\syN$ belongs to $\clSO{n}$.
We order the variables according to their subscripts: 
	i.e.,
	$x_i < x_j$ if $i<j$.
We denote by $\sbF{p}$ the subscript of the minimum variable in a Boolean polynomial $p$,
	where 
	$\sbF{p} = n+1$ if $p$ is constant.\footnote{
For instance,
	$\sbF{x_1 x_3} = 1$,
	$\sbF{x_2x_4 + x_3  + 1 } = 2$,
	and
	$\sbF{1 } = n+1$.
}	
Replacing $\plP[]{}{}$ by $\plF[]{}{}$,
	we apply \refEq{3_PL_DefPofSys} to a subset $\mathcal{P}$ in $\rgBP_n$ such that
	\envHLineCmDef[2_PL_DefFacFrm]
	{
		\plF[]{ \mathcal{P} }{}
	}
	{
		\OR{p \in \mathcal{P}} p	\lnP{=}	\Pd{p \in \mathcal{P}} (p+1) + 1		
	} 
	where $\plF[]{ \mathcal{P} }{} = 0$ if $\mathcal{P} = \phi$.
Obviously,
	\envMPd
	{
		\plP[]{ \syN }{}
	}
	{
		\plF[]{ \SetO{p_1, \ldots, p_m} }{}
	}
We define a map from the power set of $\rgBP_{n}$ to itself	
	by
	\envHLineDef[3_PL_DefMapN]
	{
		\mpN{ \mathcal{P} }
	}
	{
		\envCaseTPd{
			\Set{1}
			&
			\text{if $1 \in \mathcal{P}$}
		}{
			\mathcal{P}\setminus \Set{0}
			&
			\text{otherwise}
		}
	}	
It holds that 
	$\mpN[]{}^2 = \mpN[]{}$, 
	$\mpN{ \mathcal{P} }  \subset \mathcal{P}$		
	and									
	\envPLinePd[3_PL_EqMapNPolF]
	{
		\plF[]{ \mpN{ \mathcal{P} } }{}	\lnP{=}		\plF[]{ \mathcal{P} }{}
	}
The operations used in $\mpN[]{}$ 
	are only search of $1$ and delete of $0$.
Hence		
	the computation time of $\mpN[]{}$ is considered to be $\sbL{1}$
	by means of hashing technique (see, e.g., \cite[Section 6.4]{Knuth98} for the idea of hash).	

We are in a position to state \refThm{3_Thm1}.

\begin{theorem}	\label{3_Thm1}			
Let $\syN=\sy{p_1, \ldots, p_m}$ be a system in $\clST{n}{l}$.		

We put $\fsP{}{\tplE} = \mpN{ \Set{p_1, \ldots, p_m} }$,
	and 
	divide $\fsP{}{\tplE}$ into 
	\envHLinePd[3_Thm1_PreDef1]
	{	
		\fsP{j}{\tplE} 	
	}
	{ 
		\SetT{ p \in \fsP{}{\tplE} }{ \sbF{p} = j  } 	
	\qquad	
		(j=1, \ldots, n+1)		
	}	
For an integer $j$ from $1$ to $n+1$,
	we recursively define a family
	\envHLine[3_Thm1_FmlP]
	{
		\fmP{j} 
	}
	{	
		\SetT{  \fsP{j}{\alp_{i}\cdots\alp_{j-1}} }{ 1 \leq i \leq j-1, (\alp_{i}, \cdots, \alp_{j-1}) \in \setF[2]^{j-i} }			
	} 		
	whose elements  are subsets in $\rgBP_{j,n}$,		
	as follows.		
Firstly,
	set $\fmP{1}{} =\phi$.
Suppose $\fmP{j-1}{}$ is determined.
From the elements of $\fmP{j-1}{}$,
	we construct those of $\fmP{j}{}$ 	
	such that
	\envHLineCm[3_Thm1_DefP]
	{
		\fsP{j}{\alp_{i}\cdots\alp_{j-2}\alp_{j-1}}
	}
	{	
		\mpN{ \SetT{ p |_{x_{j-1} = \alp_{j-1} } }{ p \in \fsP{j-1}{\alp_{i}\cdots\alp_{j-2}} } }	
	}	
	where $\fsP{j-1}{\alp_{i}\cdots\alp_{j-2}} = \fsP{j-1}{\tplE}$ if $i=j-1$.	

Then, 		
	for the Boolean polynomials $\plP[]{j}{}$ in \refEq{3_PL_DefPolP} 
	with 
	$\plP[]{1}{}=\plP[]{\syN}{}$,
	we have
	\envHLinePd[3_Thm1_Main2Eq]
	{
		\plP[]{j}{} 
	}
	{
		\bkR[g]{ \Ad{\alp_{j-1} \in \setF[2]} 
			\bkR[g]{\cdots
				\bkR[g]{
					\Ad{\alp_{2} \in \setF[2]}
					\bkR[g]{
						\Ad{\alp_{1} \in \setF[2]} \plF[]{\fsP{j}{\alp_{1} \cdots \alp_{j-1}} }{} 
					}\lor \plF[]{\fsP{j}{\alp_{2} \cdots \alp_{j-1}}}{} 
				}
			\cdots  } \lor \plF[]{ \fsP{j}{\alp_{j-1}} }{}
		} 
		\lnAH
		\lor \plF[]{ \fsP{j}{\tplE} \cup \cdots \cup  \fsP{n+1}{\tplE}   }{}
	}
We also 
	have the following properties of families $\fmP{j}$.
\envItem{
\item[\prp{1}]	
	\envMPd
	{
		\nmSet{ \fmP{j} }
	}
	{
		2^{j} -2		
	}
\item[\prp{3}]
	\envMPt{\leq}
	{
		\nmSet{ \fsP{j}{\alp_{i}\cdots\alp_{j-1}} }		
	}
	{
		\nmSet{ \fsP{j-1}{\alp_{i}\cdots\alp_{j-2}} }		
	}	
	for $j \geq 2$.
\item[\prp{4}]
The computing time of \refEq{3_Thm1_DefP} 
	for all elements of $\fmP{j}$ 
	is 	
	bounded by
	\envPLinePd
	{
		\sbL[a]{ l \SmT[l]{i=1}{j-1} \Sm[l]{(\alp_i, \ldots, \alp_{j-2}) \in \setF[2]^{j-1-i} } \nmSet[a]{ \fsP{j-1}{\alp_{i}\cdots\alp_{j-2}} }  }
	}
}
\end{theorem}	

\begin{figure}[!b] \centering {\small
\caption{
The left and right trees express $\plP[]{2}{}$ and $\plP[]{3}{}$ in \refEq{3_Thm1_Main2Eq}, 
	respectively,
	where
	$\plF[]{2}{}^{\tplE}=\plF[]{ \fsP{2}{\tplE} \cup \cdots \cup \fsP{n+1}{\tplE}   }{}$
	and
	$\plF[]{3}{}^{\tplE}=\plF[]{ \fsP{3}{\tplE} \cup \cdots \cup \fsP{n+1}{\tplE}   }{}$.
}	\label{3_Fig1}	
\envPic{(400,160)
\put(0,0){
	\includegraphics[width=370pt, height=130pt, bb = 0 0 694 244]{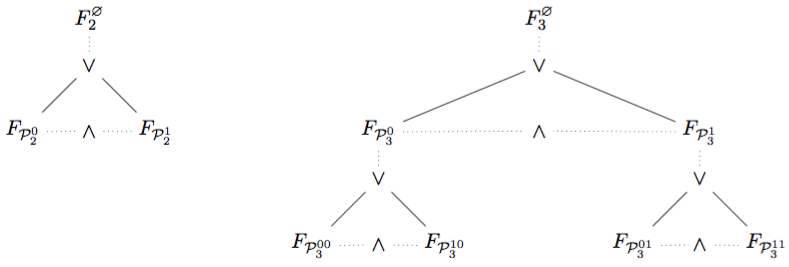}
}
}
}\end{figure}

The formula 
	\refEq{3_Thm1_Main2Eq} reads as 		
	\envHLineCm
	{
		\plP[]{1}{} 
	}
	{
		\plF[]{ \fsP{1}{\tplE} \cup \cdots \cup \fsP{n+1}{\tplE}  }{}
	}
	\vPack[20]\envHLineCm[3_PL_ExsThm2]
	{
		\plP[]{2}{} 
	}
	{
		\bkR{ 
			\plF[]{ \fsP{2}{0} }{}
			\land 
			\plF[]{ \fsP{2}{1} }{}
		} 
		\lor 
		\plF[]{ \fsP{2}{\tplE} \cup \cdots \cup \fsP{n+1}{\tplE}   }{}
	}
	\vPack[20]\envHLineCm
	{
		\plP[]{3}{} 
	}
	{
		\bkR{ 
			\bkR{
    				\bkR{ 
    					\plF[]{ \fsP{3}{00} }{}
    					\land 
    					\plF[]{ \fsP{3}{10} }{}
    				} 
    				\lor 
    				\plF[]{ \fsP{3}{0} }{}
			}
			\land 
			\bkR{
    				\bkR{ 
    					\plF[]{ \fsP{3}{01} }{}
    					\land 
    					\plF[]{ \fsP{3}{11} }{}
    				} 
    				\lor 
    				\plF[]{ \fsP{3}{1} }{}
			}
		} 
		\lor 
		\plF[]{ \fsP{3}{\tplE} \cup \cdots \cup \fsP{n+1}{\tplE}   }{}
	}
	and so on.
Both operations of conjunction and disjunction appear recursively. 		
By De Morgan's duality,
	the equations of \refEq{3_PL_ExsThm2} are equivalent to 
	\envHLineCm
	{
		\lnot \plP[]{1}{} 
	}
	{
		\lnot \plF[]{ \fsP{1}{\tplE} \cup \cdots \cup \fsP{n+1}{\tplE}  }{}
	}
	\vPack[20]\envHLineCm[3_PL_ExsThm2Dual]
	{
		\lnot \plP[]{2}{} 
	}
	{
		\bkR{ 
			\lnot \plF[]{ \fsP{2}{0} }{}
			\lor 
			\lnot \plF[]{ \fsP{2}{1} }{}
		} 
		\land
		\lnot \plF[]{ \fsP{2}{\tplE} \cup \cdots \cup \fsP{n+1}{\tplE}   }{}
	}
	\vPack[20]\envHLinePd
	{
		\lnot \plP[]{3}{} 
	}
	{
		\bkR{ 
			\bkR{
    				\bkR{ 
    					\lnot \plF[]{ \fsP{3}{00} }{}
    					\lor 
    					\lnot \plF[]{ \fsP{3}{10} }{}
    				} 
    				\land 
    				\lnot \plF[]{ \fsP{3}{0} }{}
			}
			\lor 
			\bkR{
    				\bkR{ 
    					\lnot \plF[]{ \fsP{3}{01} }{}
    					\lor 
    					\lnot \plF[]{ \fsP{3}{11} }{}
    				} 
    				\land 
    				\lnot \plF[]{ \fsP{3}{1} }{}
			}
		} 
		\land 
		\lnot \plF[]{ \fsP{3}{\tplE} \cup \cdots \cup \fsP{n+1}{\tplE}   }{}
	}
The conjunction and disjunction are replaced each other,
	and the negation is appended to each factor. 
The dual of \refEq{3_Thm1_Main2Eq} 
	is thus
	\envHLinePd[3_PL_Thm2EqDual]
	{
		\lnot \plP[]{j}{} 
	}
	{
		\bkR[g]{ \Or{\alp_{j-1} \in \setF[2]} 
			\bkR[g]{\cdots
				\bkR[g]{
					\Or{\alp_{2} \in \setF[2]}
					\bkR[g]{
						\Or{\alp_{1} \in \setF[2]} \lnot \plF[]{\fsP{j}{\alp_{1} \cdots \alp_{j-1}} }{} 
					}\land \lnot \plF[]{\fsP{j}{\alp_{2} \cdots \alp_{j-1}}}{} 
				}
			\cdots  } \land \lnot\plF[]{ \fsP{j}{\alp_{j-1}} }{}
		} 
		\lnAH
		\land \lnot \plF[]{ \fsP{j}{\tplE} \cup \cdots \cup  \fsP{n+1}{\tplE}   }{}
	}	
We can see from \refEq{3_PL_ExsThm2}
	that
	the formula \refEq{3_Thm1_Main2Eq} has expressions of binary tree as \refFig{3_Fig1},		
	in which the cases of $\plP[]{2}$ and $\plP[]{3}$ are demonstrated.
The same applies to \refEq{3_PL_Thm2EqDual}	
	with dual replacements of symbols.

Let $p$ be a Boolean polynomial in $\rgBP_n$.
We mean by $\sbE{p}$ the subscript of the maximum variable in $p$,
	where 
	$\sbE{p} = n+1$ if $p$ is constant.
For a system $\syN=\sy{p_1, \ldots, p_m}$,
	we define 
	the bandwidth 
	by
	\envHLineCmDef[3_PL_DefWthP]
	{
		\paBW_{\syN,\llaPT[]{\syN}{}}
	}
	{
		\tpMax[n]{i}{ \sbE{p_i} - \sbF{p_i} }		
	}  	
	where
	$\llaPT[]{\syN}{}$ denotes the layout of the variable order in $\syN$,
	i.e.,
	$\llaPT[]{\syN}{}$ is the map from the set of variables to $\SetO{1,\ldots,n}$		
	defined by
	$\llaPT{\syN}{x_j} = j$.
The values of $\sbF{p}$ and $\sbE{p}$
	are 
	changed 	in general 
	when variables $x_1, \ldots, x_n$ are rearranged,	
	and so
	$\paBW_{\syN,\llaPT[]{\syN}{}}$ depends on $\llaPT[]{\syN}{}$.		
Instead of $\paBW_{\syN,\llaPT[]{\syN}{}}$,
	we will use 
	$\paBW$  for short.

We will state \refThm{3_Thm2}.
The theorem comes from the fact that,
	for any $\plP[]{j}{}$ in \refEq{3_Thm1_Main2Eq},
	we can remove leaf factors as the tree depth is at most $\paBW$		
	while keeping the satisfiability.

\begin{theorem}	\label{3_Thm2}		
Let $\syN=\sy{p_1, \ldots, p_m}$ be a system in $\clST{n}{l}$,
	and
	let $\paBW$
	denote
	the bandwidth  on a variable order.	
The satisfiability of $\syN$ is decidable in time $\sbL{ 2^{\paBW} l (m + n) }$.	
\end{theorem}

We require a bit of notations on graph to state the corollaries. 	
Let $G$ be a graph,
	and
	let $V=\Fc{V}{G}$ and $E=\Fc{V}{G}$ denote the vertex and edge sets,
	respectively.
In this paper,
	we always assume 
	that 
	$G$ is simple and undirected,
	and
	that
	the vertices are ordered as $v_1, \ldots, v_n$,
	where
	$v_i < v_j$ if $i<j$.	
Let $\llaPT[]{G}{}$ denote the layout of the vertex order in $G$.
We define the bandwidth of $G$ on $\llaPT[]{G}{}$ by	
	\envHLineDefPd[3_PL_DefWthG]
	{
		\paBW_{G,\llaPT[]{G}{}}
	}
	{
		\tpMax{v_iv_j \in E}{ \abs{i-j} }
	}	
Note that
	$\paBW_{G,\llaPT[]{G}{}}=\paBW_{\syN,\llaPT[]{\syN}{}}$ 
	if 
	$\syN$ is a linear system corresponding to the adjacency matrix of $G$
	under the identifications of $x_i=v_i$.
We consider $\bkS{l}=\SetO{0, 1, \ldots, l-1}$ to be $l$ distinct colors,
	and
	we call $L$ a $\bkS{l}$-list 
	if $L \subset \bkS{l}$.

\refCor[s]{3_Cor1} and \ref{3_Cor2} are as follows.
\begin{corollary}	\label{3_Cor1}		
Let $\syN$ be a SAT problem of $n$ variables and $m$ equations,
	and
	let 
	$\paBW$ denote the bandwidth on a variable order.
\mbox{}\\{\bf (i)}
If $\syN$ is BMQ,
	we can decide the satisfiability in time $\sbL{ 2^{\paBW} (m+n) n^2 }$.
\mbox{}\\{\bf (ii)}
If $\syN$ is CNF,
	we can decide the satisfiability in time $\sbL{ 2^{\paBW} (m+n)  }$.
\end{corollary}
\begin{corollary}	\label{3_Cor2}
Let $G$ be a graph of $n$ vertices and $m$ edges,	
	and 
	let $\paBW$ denote the bandwidth on a vertex order.	
For a tuple $(L_v)_{v \in V}$ of $\bkS{l}$-lists,
	we can decide the list-colorability in time $\sbL{ (2l)^{\paBW+2} (m+ln) }$.
\end{corollary}

Our complexity results give 
	examples of FPT by bandwidth.
Our algorithms for the results are deterministic 
	as we will see in their proofs.
In practice,
	the parameters $l$ and $m$ are equal to $\sbL{n^c}$ for some constants $c$.
Particularly,
	in the list-coloring problem,
	$m$ is bounded by $\bnm{2}{n}$,
	and
	$l$ is usually considered $\sbL{1}$.

An advantage of the results is 
	that
	$f$ is either $2^{\paBW}$ or $(2l)^{\paBW+2}$
	and 
	it is not abstract,
	where
	$f$ stands for the function in the definition of FPT which is used to describe tractableness.
Therefore
	we can compare complexity results related to ours
	by ignoring
	differences in polynomial factors.

Let $\paBW \ll \Fc{g}{n}$
	mean that
	$\paBW$ is sufficiently smaller than the number $\Fc{g}{n}$
	such that 
	log values of polynomial factors in $n$ for base $2$ 
	have 
	no influence.
In \cite{LPTWY17},
	a randomized algorithm for BMQ-SAT	
	is presented,
	whose computation time is bounded by $\sbLp{2^{0.8765n}}$.
Our algorithm for BMQ
	is faster 
	when $\paBW \ll 0.8765n$.
In \cite{Hertli14A},
	randomized algorithms for $3$-CNF-SAT and $4$-CNF-SAT are presented,
	whose computation times are bounded by
	$\sbLp{1.30704^n}$ and $\sbLp{1.46899^n}$,
	respectively. 		
Our algorithm for CNF
	is faster 
	when $\paBW \ll 0.38630 n$ and $\paBW \ll 0.55482 n$,		
	because
	$1.30704 \fallingdotseq 2^{0.38630}$
	and
	$1.46899 \fallingdotseq 2^{0.55482}$,
	respectively.
In \cite{BHK09},	
	a combinatorial algorithm for list-coloring is presented,	
	whose computation time is bounded by $\sbLp{2^{n}}$.
Our algorithm for list-coloring
	is faster 
	when $\paBW \ll n$,		
	where 
	we assume
	that
	the number $l$ of colors 	
	is 
	constant and independent to $n$.

We end the section with an additional comparison.	
In \cite{JS97} (see also \cite{FFLRSST11}),	
	a dynamic programming algorithm for list-coloring is presented,
	whose computation time is bounded by $\sbL{n^{\paTW+2}}$,	
	where $\paTW$ is treewidth.
Let $\sbLs[]{}$
	denote the little-o notation.
Since
	$m=\sbL{n^2}$,
	our algorithm
	is faster 
	when 
	\envMPd{
		(2l)^{\paBW} 
	}
	{ 
		\sbLs{ n^t }
	}


\section{Proof of \refThm{3_Thm1}}	\label{sectFour}
We require \refLem[s]{4_Lem1} and \ref{4_Lem2} to prove \refThm{3_Thm1}.

\begin{lemma}	\label{4_Lem1}
For positive integers $h$ and $j$ with $h \leq j$,
	we have
	\envHLinePd[4_Lem1_Eq1]
	{
		\SmT{i=h}{j-1}  \nmSet[a]{ \setF[2]^{j-i} }
	}
	{
		2^{j+1-h} - 2
	}
\end{lemma}
\begin{lemma}	\label{4_Lem2}
Let $p$ be a Boolean polynomial in $\spCLT{n}{l}$, 	
	let $x_h$ be a variable,
	and
	let $\alp$ be a value in $\setF[2]$.
The computing time of $p |_{x_h = \alp}$ is in $\sbL{l}$.	
\end{lemma}

We will first prove \refThm{3_Thm1} dividing into two parts: 	
	one is devoted to the formula \refEq{3_Thm1_Main2Eq}
	and
	the other is devoted to the properties \prp{1}, \prp{3} and \prp{4}.	
We will then prove \refLem[s]{4_Lem1} and \ref{4_Lem2}.

\envProof[the formula \refEq{3_Thm1_Main2Eq}]{	
We will use induction on $j$. 
The case of $P_1$ 
	is
	obvious,	
	because
	\envPLinePd
	{
		P_1 
	\lnP{=}	
		\plP[]{\sy{p_1, \ldots, p_m}}{}
	\lnP{ \osTx{=}{\refEq{3_PL_DefPofSys} \atop \refEq{2_PL_DefFacFrm}} }			
		\plF[]{ \Set{p_1, \ldots, p_m} }{}		
	\lnP{\osTx{=}{\refEq{3_PL_EqMapNPolF}  }}
		\plF[]{ \mpN{ \Set{p_1, \ldots, p_m} } }{}
	\lnP{ \osTx{=}{ \refEq{3_Thm1_PreDef1}} }
		\plF[]{ \fsP{1}{\tplE} \cup \cdots \cup \fsP{n+1}{\tplE} }{}
	}
Let $j > 1$,
	and
	suppose that 
	\refEq{3_Thm1_Main2Eq} is true in the case of $P_{j-1}$.
Let $\alp_{j-1} \in \setF[2]$.	
We have
	\envHLine
	{
		p |_{x_{j-1} = \alp_{j-1} } 
	}
	{ 
		p 
	} 
	for a Boolean polynomial $p$ with $\sbF{p} > j-1$,
	and so 
	\envHLineCm
	{
		\SetT{  p |_{x_{j-1} = \alp_{j-1} } }{ p \in \fsP{j-1}{\tplE} \cup \cdots \cup \fsP{n+1}{\tplE} }
	}
	{
		\SetT{ p |_{x_{j-1} = \alp_{j-1} } }{ p \in \fsP{j-1}{ \tplE }  } 
		\lnP{\sqcup}
		\bkR{ \fsP{j}{\tplE} \cup \cdots \cup \fsP{n+1}{\tplE} }
	}
	where
	$\sqcup$ means the disjoint union.
Hence,	
	by \refEq{2_PL_DefFacFrm},
	\envHLineFiCm
	{
		\plF[]{ \fsP{j-1}{\tplE} \cup \cdots \cup \fsP{n+1}{\tplE}  }{} \mVert[B]_{x_{j-1}=\alp_{j-1}}
	}
	{ 
		\OR{p \in \fsP{j-1}{\tplE} \cup \cdots \cup \fsP{n+1}{\tplE} } p	\mVert_{x_{j-1}=\alp_{j-1}}
	}
	{
		\OR{p' \in \SetT{  p |_{x_{j-1} = \alp_{j-1} } }{ p \in \fsP{j-1}{\tplE} \cup \cdots \cup \fsP{n+1}{\tplE} } } p'	
	}
	{
		\bkR[a]{ \OR{p' \in \SetT{ p |_{x_{j-1} = \alp_{j-1} } }{ p \in \fsP{j-1}{ \tplE }  }   } p' }
		\lor
		\bkR[a]{ \OR{p' \in \fsP{j}{\tplE} \cup \cdots \cup \fsP{n+1}{\tplE}  } p'	}
	}
	{\rule{0pt}{14pt}
		\plF[]{ \SetT{ p |_{x_{j-1} = \alp_{j-1} } }{ p \in \fsP{j-1}{ \tplE }  }  }{} \lor \plF[]{ \fsP{j}{\tplE} \cup \cdots \cup \fsP{n+1}{\tplE} }{}
	}
	which,
	together with \refEq{3_PL_EqMapNPolF}  and \refEq{3_Thm1_DefP}, 
	yields
	\envHLinePd[4_Pr_3_Thm1_Eq1Help]	
	{
		\plF[]{ \fsP{j-1}{\tplE} \cup \cdots \cup \fsP{n+1}{\tplE}  }{} \mVert[B]_{x_{j-1}=\alp_{j-1}}
	}
	{	
		\plF[]{ \fsP{j}{\alp_{j-1} } }{} \lor \plF[]{ \fsP{j}{\tplE} \cup \cdots \cup \fsP{n+1}{\tplE} }{}
	}
For an element $\fsP{j-1}{\alp_{i}\cdots\alp_{j-2}}$ in $\fmP{j-1}$,
	combining \refEq{2_PL_DefFacFrm}, \refEq{3_PL_EqMapNPolF} and \refEq{3_Thm1_DefP}
	also yields
	\envHLinePd[4_Pr_3_Thm1_Eq2Help]	
	{
		\plF[]{ \fsP{j-1}{\alp_{i}\cdots\alp_{j-2}} }{} \mVert[B]_{x_{j-1}=\alp_{j-1}}
	}
	{
		\plF[]{ \fsP{j}{\alp_{i}\cdots\alp_{j-2} \alp_{j-1} } }{}
	}
By the induction hypothesis,
	$\plP[]{j-1}{}$  satisfies \refEq{3_Thm1_Main2Eq}.
Therefore,
	by \refEq{4_Pr_3_Thm1_Eq1Help} and \refEq{4_Pr_3_Thm1_Eq2Help},
	\envLineThCm[4_Pr_3_Thm1_Eq2]
	{
		\plP[]{j-1}{} |_{x_{j-1}=\alp_{j-1}} 
	}
	{
        		\bkR[g]{ \Ad{\alp_{j-2} \in \setF[2]} 
                    	\bkR[g]{\cdots
        			\bkR[g]{
        				\Ad{\alp_{2} \in \setF[2]}
                    			\bkR[g]{
                    				\Ad{\alp_{1} \in \setF[2]} \plF[]{\fsP{j-1}{\alp_{1} \cdots \alp_{j-2}  }}{} 
                    			} \lor \plF[]{\fsP{j-1}{ \alp_{2} \cdots \alp_{j-2}  }}{} 
        			}
                    	\cdots  } \lor \plF[]{ \fsP{j-1}{ \alp_{j-2}  } }{}
        		} \mVert[g]_{x_{j-1}=\alp_{j-1}} 	
		\lnAH
		\lor \plF[]{ \fsP{j-1}{\tplE} \cup \cdots \cup \fsP{n+1}{\tplE}  }{} \mVert[B]_{x_{j-1}=\alp_{j-1}}
	}
	{
		G_{j, \alp_{j-1}}
		\lor	
		\plF[]{ \fsP{j}{\alp_{j-1} } }{} 
		\lor 
		\plF[]{ \fsP{j}{\tplE} \cup \cdots \cup \fsP{n+1}{\tplE}  }{} 
	}
	where 
	\envLinePd
	{
		G_{j, \alp_{j-1}}
	}
	{
		\Ad{\alp_{j-2} \in \setF[2]} 
                	\bkR[g]{\cdots
			\bkR[g]{
				\Ad{\alp_{2} \in \setF[2]}
                			\bkR[g]{
                				\Ad{\alp_{1} \in \setF[2]} \plF[]{\fsP{j}{\alp_{1} \cdots \alp_{j-2} \alp_{j-1} }}{} 
                			} \lor \plF[]{\fsP{j}{\alp_{2} \cdots \alp_{j-2} \alp_{j-1} }}{} 
			}
                	\cdots  } \lor \plF[]{ \fsP{j}{\alp_{j-2} \alp_{j-1} } }{}
	}
Because of \refEq{3_PL_DefPolP} and \refEq{4_Pr_3_Thm1_Eq2},
	$\plP[]{j}{}$ is expressed as
	\envHLinePd
	{
		\plP[]{j}{} 
	}
	{
		\bkR{ G_{j,0} \lor \plF[]{ \fsP{j}{ 0 } }{} \lor \plF[]{ \fsP{j}{\tplE} \cup \cdots \cup \fsP{n+1}{\tplE} }{} }
		\land
		\bkR{ G_{j,1} \lor \plF[]{ \fsP{j}{ 1 } }{} \lor \plF[]{ \fsP{j}{\tplE} \cup \cdots \cup \fsP{n+1}{\tplE} }{} }
	}
Using the distributivity of $\lor$ over $\land$,	
	we obtain
	\envHLineThCm[4_Pr_3_Thm1_Eq3]
	{
		\plP[]{j}{} 
	}
	{
		\bkR{
			\bkR{ G_{j,0} \lor \plF[]{ \fsP{j}{ 0 } }{} }
			\land
			\bkR{ G_{j,1} \lor \plF[]{ \fsP{j}{ 1 } }{} }
		}
		\lor \plF[]{ \fsP{j}{\tplE} \cup \cdots \cup \fsP{n+1}{\tplE} }{} 
	}
	{
		\bkR[g]{ \Ad{ \alp_{j-1} \in \setF[2]} \bkR{ G_{j,\alp_{j-1}} \lor \plF[]{ \fsP{j}{ \alp_{j-1} } }{} } }	
		\lor	
		\plF[]{ \fsP{j}{\tplE} \cup \cdots \cup \fsP{n+1}{\tplE} }{} 
	}
	which shows that
	\refEq{3_Thm1_Main2Eq} is true in the case of $P_{j}$.
}
\envProof[the properties \prp{1}, \prp{3} and \prp{4}]{	
The property \prp{1} immediately follows from \refEq{3_Thm1_FmlP} and \refEq{4_Lem1_Eq1} with $h=1$.
 	
For an element $\fsP{j}{\alp_{i}\cdots\alp_{j-1}}$ in $\fmP{j}$ for $j \geq 2$,	
	we have		
	\envHLineCmPt{=}	
	{
		\nmSet[a]{ \fsP{j}{\alp_{i}\cdots\alp_{j-1}} }
	}
	{
		\nmSet[a]{ \mpN{ \SetT{ p |_{x_{j-1} = \alp_{j-1} } }{ p \in \fsP{j-1}{\alp_{i}\cdots\alp_{j-2}} } } }
	\lnAHP{\leq}	
		\nmSet[a]{ \SetT{ p |_{x_{j-1} = \alp_{j-1} } }{ p \in \fsP{j-1}{\alp_{i}\cdots\alp_{j-2}} } }
	\lnAHP{\leq}
		\nmSet[a]{\fsP{j-1}{\alp_{i}\cdots\alp_{j-2}} }
	}
	which implies \prp{3}.

For an element $\fsP{j}{\alp_{i}\cdots\alp_{j-1}}$,		
	we can calculate 
	\refEq{3_Thm1_DefP} 
	in time
	\envMO
	{
		\sbL{l \nmSet{ \fsP{j-1}{\alp_{i}\cdots\alp_{j-2}} } }
	}
	by \refLem{4_Lem2},
	where remember that
	$\mpN[]{}$ is a constant cost map 
	by hashing technique.		
Thus,
	by \refEq{3_Thm1_FmlP},
	the computing time of \refEq{3_Thm1_DefP} for all elements of $\fmP{j}$ 
	is bounded by
	\envHLineCm
	{
		\sbL[a]{ l \SmT[l]{i=1}{j-1} \Sm[l]{(\alp_i, \ldots, \alp_{j-1}) \in \setF[2]^{j-i} } \nmSet[a]{ \fsP{j-1}{\alp_{i}\cdots\alp_{j-2}} }  }
	}
	{
		\sbL[a]{ 2 l \SmT[l]{i=1}{j-1} \Sm[l]{(\alp_i, \ldots, \alp_{j-2}) \in \setF[2]^{j-1-i} } \nmSet[a]{ \fsP{j-1}{\alp_{i}\cdots\alp_{j-2}} }  }
	}
	which proves \prp{4}. 		
}

We will prove \refLem{4_Lem1}.		

\envProof[\refLem{4_Lem1}]{
We have
	\envOFLineFCm
	{
		\SmT{i=h}{j-1}  \nmSet[a]{ \setF[2]^{j-i} }
	}
	{
		\SmT{i=h}{j-1} 2^{j-i}
	}
	{
		2 \bkR{ 1 + 2 + \cdots + 2^{j-1-h} }
	}
	{
		2 \bkR{ 2^{j-h} -1 }
	}
	which implies \refEq{4_Lem1_Eq1}.
}

Let $\stL_n$ be the set 
	consisting of the literals 
	and
	the constant $1$.
We define 
	a map $\mpE[]{}$ from $\spCLO{n}$ to the power set of $\stL_n$ by 
	\envHLineDef[4.2_PL_DefEmbMap]
	{
		\mpE{c}
	}
	{
		\envCaseThCm{
			\phi
			&
			\text{if $c=0$}
		}{
			\SetO{1}
			&
			\text{if $c=1$}
		}{
			\SetO{ l_1, \ldots, l_k }
			&
			\text{if $c = l_1 \cdots l_k$}
		}
	}
	where $c \in \spCLO{n}$ and $l_i$ are literals. 	
This map is well-defined and injective by \refEq{2_PL_AsmCls}.	

We will prove \refLem{4_Lem2}.

\envProof[\refLem{4_Lem2}]{
We may show
	that
	the computing time of $c |_{x_h = \alp}$ is in $\sbL{1}$
	for a non-constant clause $c$.
Let $k$ denote the degree of $c$.	
There exist
	$k$ variables $x_{h_i}$
	and
	$k$ values $\beta_i$ in $\setF[2]$ such that 
	\envHLineCm[Pr_4_Lem2_Eq1]
	{
		c
	}
	{
		\PdT{i=1}{k}(x_{h_i}+\beta_i)
	}
	where
	literals $x_{h_i}+\beta_i$ satisfy \refEq{2_PL_AsmCls}.	
Put $\mathcal{X} = \SetO{x_{h_1}, \ldots, x_{h_k}}$.
If $x_h \notin \mathcal{X}$,
	\envHLinePd[Pr_4_Lem2_Eq2]
	{
		c |_{x_h = \alp}
	}
	{ 
		c
	}
If $x_h \in \mathcal{X}$,
	\envHLineF[Pr_4_Lem2_Eq3]
	{
		c |_{x_{h} = \alp }
	}
	{
		(x_h + \beta_h) |_{x_h = \alp} \tpPdT{i=1}{h_i \neq h}{k}(x_{h_i}+\beta_i)
	}
	{
		(\alp + \beta_h) \tpPdT{i=1}{h_i \neq h}{k}(x_{h_i}+\beta_i)
	}
	{
		\envCaseTPd{
			0
			&
			\text{if $\alp = \beta_h$}
		}{
			\tpPdT[l]{i=1}{h_i \neq h}{k}(x_{h_i}+\beta_i)
			&
			\text{if $\alp \neq \beta_h$}
		}
	}
We can see 
	from \refEq{Pr_4_Lem2_Eq1}, \refEq{Pr_4_Lem2_Eq2} and \refEq{Pr_4_Lem2_Eq3} 
	that	
	evaluating $\mpE{ c |_{x_h = s + \alp} }$ from $\mpE{c}$ 
	is 
	implemented by the following process:
\envItemEm{
\item
Set $\ome = \mpE{c}$.
\item
Return $\ome$ 
	if 
	$\ome = \phi$
	or
	$1 \in \ome$.\footnote{
Note that
	$\ome = \phi$ if and only if $c=0$.
Also 
	note that
	$1 \in \ome$ if and only if $c=1$;
	the reason 
	is because 	 
	$\mpE{c} = \SetO{1}$ if $c=1$
	and
	$1 \notin \mpE{c}$ otherwise.	
Hence 
	the return condition of step \theenumi\mbox{} is equivalent to $c \in \setF[2]$.
}
\item
Search $x_h$ and $x_h + 1$ from $\ome$.
Return $\ome$ if not exist.
\item
\newcounter{tmpNum}	\setcounter{tmpNum}{\theenumi} 	\addtocounter{tmpNum}{-1}
Set $l_h = x_h + \beta_h = \text{(the literal searched in step \arabic{tmpNum})}$.\footnote{		
Note that 
	$\nmSet{ \SetO{x_h, x_h+1} \cap \mpE{c} } \leq 1$ by \refEq{2_PL_AsmCls},
	and
	$l_h$ in step \theenumi\mbox{} is uniquely determined.
}
\item
Delete $l_h$ from $\ome$.
\item
Return $\phi$ if $\alp=\beta_h$,	
	otherwise 
	return $\ome$. 	
}
The operations 
	used in the process which are not elemental 
	are
	search and delete. 	
By hashing technique,
	costs of these operations are constants.
Thus
	the computation time of the process  is bounded by $\sbL{1}$.
Since
	$\mpE[]{}$ 
	is 
	an embedding of $\spCLO{n}$,
	the process implies 
	that
	the time of computing $c |_{x_h = \alp}$ 
	is bounded by $\sbL{1}$,
	and
	we compete the proof.
}

\section{Proof of \refThm{3_Thm2}}	\label{sectFive}
We will require \refProp{5_Prp1} to prove \refThm{3_Thm2}.
The proposition
	is a refinement of \refThm{3_Thm1},
	which has 
	an additional condition of the bandwidth.	%

\begin{proposition}	\label{5_Prp1}	
Let $\syN=\sy{p_1, \ldots, p_m}$ be a system in $\clST{n}{l}$ 
	with the bandwidth $\paBW$ on a variable order. 
Let $\fsP{j}{\tplE}$ denote 
	the subsets defined in \refEq{3_Thm1_PreDef1},
	and
	let $\fmP{\paBW+1}$  
	denote
	the $(\paBW+1)$-th family determined in \refEq{3_Thm1_FmlP} and \refEq{3_Thm1_DefP}.

We put
	\envM
	{
		I_j 
	}
	{ 
		j-\paBW	
	}	
	for 
	$j \geq \paBW + 1$.		
For an integer $j$ from $\paBW+1$ to $n+1$,
	we recursively define a family
	\envHLine[5_Prp1_FmlQ]
	{
		\fmQ{j} 
	}
	{	
		\SetT{  \fsQ{j}{\alp_{i}\cdots\alp_{j-1}} }{ I_j \leq i \leq j-1, (\alp_{i}, \cdots, \alp_{j-1}) \in \setF[2]^{j-i} }		
	} 		
	whose elements  are subsets in $\rgBP_{j,n}$,		
	as follows.		
Firstly,
	set
	\envMPd
	{
		\fmQ{\paBW+1} 
	}
	{ 
		\fmP{\paBW+1}
	} 
Suppose 
	that
	$\fmQ{j-1}{}$ is determined.	
From the elements of $\fmQ{j-1}{}$,
	we construct temporal elements such that		
	\envHLineCm[5_Prp1_DefQm]
	{
		\fsQ[m]{j}{\alp_{i}\cdots\alp_{j-2}\alp_{j-1}}
	}
	{
		\mpN{ \SetT{ p |_{x_{j-1} = \alp_{j-1} } }{ p \in \fsQ{j-1}{\alp_{i}\cdots\alp_{j-2}} } }
	}	
	where
	$\fsQ{j-1}{\alp_{i}\cdots\alp_{j-2}} = \fsP{j-1}{\tplE}$ if $i=j-1$.
Then
	we define the elements of $\fmQ{j}{}$ by\,\footnote{
We can replace ``$>0$'' with ``$=1$'' in \refEq{5_Prp1_DefQ},
	because 
	$\fsQ[m]{j}{\alp_{I_{j-1}} \alp_{I_{j}}\cdots\alp_{j-1}} \subset \SetO{1}$ by \refEq{5_Pr_Prp1_Eq8}
	and
	the number of elements in $\fsQ[m]{j}{\alp_{I_{j-1}} \alp_{I_{j}}\cdots\alp_{j-1}}$ 
	is one or zero.
}
	\envHLineDef[5_Prp1_DefQ]
	{
		\fsQ{j}{\alp_{i}\cdots\alp_{j-1}}
	}
	{
		\envCaseTPd	
		{	
			\SetO{1}
			&	
			\text{if $i=I_j$ 	
				and 
				$\Pd[l]{\alp_{I_{j-1}} \in \setF[2]} \nmSet[a]{\fsQ[m]{j}{\alp_{I_{j-1}} \alp_{I_{j}}\cdots\alp_{j-1}} } > 0$} 	
		}		
		{	
			\fsQ[m]{j}{\alp_{i}\cdots\alp_{j-1}}
			&		
			\text{if $i>I_j$ 		
				or 
				$\Pd[l]{\alp_{I_{j-1}} \in \setF[2]} \nmSet[a]{\fsQ[m]{j}{\alp_{I_{j-1}} \alp_{I_{j}}\cdots\alp_{j-1}} } = 0$}	
		}
	}		

Let
	$j \geq \paBW+1$.
Then, 
	for the Boolean polynomials $\plP[]{j}{}$	
	in \refEq{3_PL_DefPolP} with $\plP[]{1}{}=\plP[]{\syN}{}$,
	we have
	\envHLinePd[5_Prp1_Main2Eq]
	{
		\plP[]{j}{} 
	}
	{
		\bkR[g]{ \Ad{\alp_{j-1} \in \setF[2]} 
			\bkR[g]{\cdots
				\bkR[g]{
					\Ad{\alp_{I_{j}+1} \in \setF[2]}
					\bkR[g]{
						\Ad{\alp_{I_{j}} \in \setF[2]} \plF[]{\fsQ{j}{\alp_{I_{j}} \cdots \alp_{j-1}} }{} 
					}\lor \plF[]{\fsQ{j}{\alp_{I_{j}+1} \cdots \alp_{j-1}}}{} 
				}
			\cdots  } \lor \plF[]{ \fsQ{j}{\alp_{j-1}} }{}
		} 
		\lnAH
		\lor \plF[]{ \fsP{j}{\tplE} \cup \cdots \cup \fsP{n+1}{\tplE} }{}
	}
We also 
	have the following properties of families $\fmQ{j}$.
\envItem{
\item[\prpS{1}]	
	\envMPd	
	{
		\nmSet{ \fmQ{j} }
	}
	{
		2^{ \paBW+1 } -2			
	}
\item[\prpS{3}]
	If $i>I_j$,
	\envMPdPt{\leq}
	{
		\nmSet{ \fsQ{j}{\alp_{i}\cdots\alp_{j-1}} }		
	}
	{		
		\nmSet{ \fsQ{j-1}{\alp_{i}\cdots\alp_{j-2}} }		
	}
	If $i=I_j$,
	\envMPdPt{\leq}
	{
		\nmSet{ \fsQ{j}{\alp_{I_j}\cdots\alp_{j-1}} }
	}
	{
		2
	}
\item[\prpS{4}]
The computing time of \refEq{5_Prp1_DefQm} and \refEq{5_Prp1_DefQ} for all elements of $\fmQ{j}$ 	
	is bounded by
	\envPLinePd
	{
		\sbL[a]{ 2^{\paBW} 
		+  
		l \SmT[l]{i=I_{j-1}}{j-1} \Sm[l]{ (\alp_i, \ldots, \alp_{j-2}) \in \setF[2]^{j-1-i} } \nmSet[a]{ \fsQ{j-1}{\alp_{i}\cdots\alp_{j-2}} }  }
	}
}
\end{proposition}

We will prove \refThm{3_Thm2}.
Then we will prove \refProp{5_Prp1}.

\envProof[\refThm{3_Thm2}]{	
Let $\fmP{1}, \ldots, \fmP{\paBW+1}$  
	denote
	the first $(\paBW+1)$ families determined in \refEq{3_Thm1_FmlP} and \refEq{3_Thm1_DefP}.
Let $j \in \SetO{2, \ldots, \paBW+1}$. 		
Using \prp{3} in \refThm{3_Thm1} repeatedly,
	we obtain 
	\envPLine
	{
		\nmSet[a]{ \fsP{j-1}{\alp_{i}\cdots\alp_{j-2}} }
	\lnP{\leq}
		\nmSet[a]{ \fsP{j-2}{\alp_{i}\cdots\alp_{j-3}} }
	\lnP{\leq}
		\cdots
	\lnP{\leq}
		\nmSet[a]{ \fsP{i+1}{\alp_{i}} }
	\lnP{\leq}
		\nmSet[a]{ \fsP{i}{\tplE} }
	}
	for
	any $\fsP{j-1}{\alp_{i}\cdots\alp_{j-2}} \in \fmP{j-1}$.
Since
	\envM{
		\nmSet{\setF[2]^{h}}
	}{
		2^{h}
	}
	for $h\geq 0$,
	we have
	\envHLineCmPt[5_Pr_3_Thm2_Eq1]{ \leq }	
	{
		\SmT{i=1}{j-1} \Sm[l]{(\alp_i, \ldots, \alp_{j-2}) \in \setF[2]^{j-1-i} } \nmSet[a]{ \fsP{j-1}{\alp_{i}\cdots\alp_{j-2}} }
	}
	{
		\SmT{i=1}{j-1} 2^{j-1-i} \nmSet[a]{ \fsP{i}{\tplE} }
	}	
	which, 
	together with \prp{4} in \refThm{3_Thm1},
	shows that
	the total time to calculate the families $\fmP{2}, \ldots, \fmP{\paBW+1}$		
	is bounded by
	\envPLinePd[5_Pr_3_Thm2_Eq3]
	{
		\sbL[a]{ l \SmT{j=2}{\paBW+1}\SmT{i=1}{j-1} 2^{j-1-i} \nmSet[a]{ \fsP{i}{\tplE} } }	
	}
Note that
	$\fmP{1} = \phi$ by the initial condition
	and 
	no calculation is required for $\fmP{1}$.

Let 
	$\fmQ{\paBW+2}, \ldots, \fmQ{n+1}$ 
	denote 
	the families 
	determined in \refEq{5_Prp1_FmlQ}, \refEq{5_Prp1_DefQm} and \refEq{5_Prp1_DefQ} 
	with the initial condition $\fmQ{\paBW+1} = \fmP{\paBW+1}$.
Let $j \in \SetO{\paBW +2, \ldots, n+1}$.	
Similarly to \refEq{5_Pr_3_Thm2_Eq1},		
	it follows 
	from	$I_{j-1} = I_{j}  -1 = j-1-\paBW$	and \prpS{3} in \refPrp{5_Prp1} 	
	that
	\envLineCmPt{ = }	
	{
		\SmT{i=I_{j-1}}{j-1} \Sm[l]{(\alp_i, \ldots, \alp_{j-2}) \in \setF[2]^{j-1-i} } \nmSet[a]{ \fsQ{j-1}{\alp_{i}\cdots\alp_{j-2}} }
	}
	{
		\Sm{ (\alp_{I_{j-1}}, \ldots, \alp_{j-2}) \in \setF[2]^{\paBW} } \nmSet[a]{ \fsQ{j-1}{\alp_{I_{j-1}}\cdots\alp_{j-2}} }
		+
		\SmT{i=I_j}{j-1} \Sm[l]{(\alp_i, \ldots, \alp_{j-2}) \in \setF[2]^{j-1-i} } \nmSet[a]{ \fsQ{j-1}{\alp_{i}\cdots\alp_{j-2}} }
	\lnAHP{\leq}
		2\cdot 2^{\paBW} + \SmT[l]{i=j-\paBW}{j-1} 2^{j-1-i} \nmSet[a]{ \fsP{i}{\tplE} }
	}	
	which, 
	together with \prpS{4} in \refPrp{5_Prp1}, 	
	shows
	that
	the total time to calculate the families $\fmQ{\paBW+2}, \ldots, \fmQ{n+1}$
	is bounded by
	\envPLinePd[5_Pr_3_Thm2_Eq4]
	{
		\sbL[a]{ 
			l n  2^{\paBW} 
			+ 
			l \SmT{j=\paBW+2}{n+1}\SmT{i=j-\paBW}{j-1} 2^{j-1-i} \nmSet[a]{ \fsP{i}{\tplE} } 
		}	
	}

We define subsets in $\setZ^2$ as follows:
	\envPLine
	{\begin{array}{ccccrclrcl}
		L_1 &:=& \{ (i,j)\in\setZ^2 & | & 2  	& \leq j \leq& \paBW+1, &1 				&\leq i \leq  &j-1 	\}, 	\\
		L_2 &:=& \{ (i,j)\in\setZ^2 & | & \paBW+2  	& \leq j \leq& n+1, &j-\paBW 		&\leq i \leq  &j-1 	\},	\rule{0pt}{15pt}\\
		L_3 &:=& \{ (i,j)\in\setZ^2 & | & n+2  	& \leq j \leq& n+\paBW+1, \,\, &j-\paBW	&\leq i \leq  &n+1  		\}.	\rule{0pt}{15pt}
	\end{array}}
Obviously,
	$L_a \cap L_b = \phi$ for $a \neq b$.
Switching the roles of $i$-axis and $j$-axis,
	we obtain
	\envPLine
	{\begin{array}{ccccrclrcl}
		L_1 &=& \{ (i,j)& \in\setZ^2 \,\, | & 1  	& \leq i \leq& \paBW, 			& i+1 	&\leq j \leq  &\paBW+1 	\}, 	\\
		L_2 &=& \{ (i,j)& \in\setZ^2 \,\, | & 2  	& \leq i \leq& \paBW, 			&\paBW+2 	&\leq j \leq  &i+\paBW 	\}	\rule{0pt}{15pt}\\
		       &  &\cup \{ (i,j)& \in\setZ^2 \,\, | & \paBW+1  	& \leq i \leq& n+1-\paBW, \,\,	&i+1 	&\leq j \leq  &i+\paBW 	\}	\rule{0pt}{12pt}\\
		       &  &\cup \{ (i,j)& \in\setZ^2 \,\, | & n+2-\paBW  	& \leq i \leq& n, &i+1 	&\leq j \leq  &n+1 	\},	\rule{0pt}{12pt}\\
		L_3 &=& \{ (i,j)& \in\setZ^2 \,\, | & n+2-\paBW  	& \leq i \leq& n+1, &n+2	&\leq j \leq  &i+\paBW  		\}.	\rule{0pt}{15pt}
	\end{array}}
Hence,
	\envHLineCm
	{
		L_1 \cup L_2 \cup L_3
	}
	{
		\SetT{ (i,j)\in\setZ^2 }{ 1 \leq i \leq n+1, \, i+1 \leq j \leq i+\paBW}
	}
	and
	\envHLinePd
	{
		\SmT{j=2}{\paBW+1}\SmT{i=1}{j-1} 2^{j-1-i} \nmSet[a]{ \fsP{i}{\tplE} } 
		+
		\SmT{j=\paBW+2}{n+1}\SmT{i=j-\paBW}{j-1} 2^{j-1-i} \nmSet[a]{ \fsP{i}{\tplE} }
	}
	{
		\Sm{ (i,j) \in L_1 \cup L_2 } 2^{j-1-i} \nmSet[a]{ \fsP{i}{\tplE} }		
	\lnAHP{\leq}	
		\Sm{ (i,j) \in L_1 \cup L_2 \cup L_3} 2^{j-1-i} \nmSet[a]{ \fsP{i}{\tplE} }
	\lnAHP{=}	
		\SmT{i=1}{n+1}  \nmSet[a]{ \fsP{i}{\tplE} }  \SmT{j=i+1}{i+\paBW} 2^{j-1-i}		
	}
Since
	\envPLineCm
	{&&
		\SmT[l]{i=1}{n+1} \nmSet{ \fsP{i}{\tplE} } 
	\lnP{=}
		\nmSet{\mpN{ \Set{p_1, \ldots, p_m} }} 
	\lnP{\leq} 
		m
	,\qquad
		\SmT[l]{j=i+1}{i+\paBW} 2^{j-1-i} 
	\lnP{=} 
		\SmT[l]{j=1}{\paBW} 2^{j-1} 
	\lnP{<} 
		2^{\paBW}
	}
	we have
	\envHLinePdPt[5_Pr_3_Thm2_Eq5]{<}
	{
		\SmT{j=2}{\paBW+1}\SmT{i=1}{j-1} 2^{j-1-i} \nmSet[a]{ \fsP{i}{\tplE} } 
		+
		\SmT{j=\paBW+2}{n+1}\SmT{i=j-\paBW}{j-1} 2^{j-1-i} \nmSet[a]{ \fsP{i}{\tplE} }
	}
	{
		m 2^{\paBW}
	}
Therefore,
	we see from 
	\refEq{5_Pr_3_Thm2_Eq3}, \refEq{5_Pr_3_Thm2_Eq4}, and \refEq{5_Pr_3_Thm2_Eq5}
	that
	the whole time to calculate all families $\fmP{i}$ and $\fmQ{j}$
	is bounded by
	\envPLinePd[5_Pr_3_Thm2_Eq6]
	{
		\sbL{ l (m + n ) 2^{\paBW} }
	}
	
It is required to compute \refEq{3_Thm1_PreDef1} 
	for starting the above procedure to calculate all the families;
	this
	costs in $\sbL{lm}$
	since 
	\refEq{3_Thm1_PreDef1} is done by dividing $m$ polynomials consisting of $l$ clauses 
	into $(n+1)$ sets. 
The solvability of $\syN$ is equivalent to $\plP[]{n+1}{}=0$,	
	and
	it is also required to confirm whether $\plP[]{n+1}{}$ 
	is zero or not 
	for closing the procedure;
	this
	costs in $\sbL{2^{\paBW}}$,
	since,
	by \prpS{1},
	the number of factors in the right-hand side of \refEq{5_Prp1_Main2Eq} for $j=n+1$ 
	is 
	less than $2\cdot 2^{\paBW}$.\footnote{
Note that
	the factors belong to $\setF[2]$ because $\plP[]{n+1}{} \in \setF[2]$,
	and
	that
	binary operations on $\setF[2]$ cost in $\sbL{1}$.
}
Both computation times for starting and closing
	are bounded by \refEq{5_Pr_3_Thm2_Eq6},
	and
	we prove \refThm{3_Thm2}.
}

We require the following lemmas to show \refPrp{5_Prp1}.		

\begin{lemma}	\label{5_Lem1}
Let $p$ be a Boolean polynomial,
	and
	let $\paBW$ be a positive integer
	such that
	$\sbE{p} - \sbF{p} \leq \paBW$.		
Put $i=\sbF{p}$,
	and
	let $j$ be an integer with $i < j \leq n+1$.	
Then	 we have
	\envHLinePt[5_Lem1_Eq1]{\in}
	{
		p |_{x_{i} = \alp_{i} } |_{x_{i+1} = \alp_{i+1} } \cdots  |_{x_{j-1} = \alp_{j-1} }	
	}
	{
		\rgBP_{j, i+\paBW }
	}	
	for values $\alp_{i}, \alp_{i+1}, \ldots, \alp_{j-1}$ in $\setF[2]$.		
\end{lemma}
\begin{lemma}	\label{5_Lem2}
Let $\mathcal{P}$ be a subset in $\rgBP_{n}$.
Let $i,j$ be integers with $1 \leq i < j \leq n+1$,
	and 
	let 
	$\alp_i, \ldots, \alp_{j-1}$ 
	be values in $\setF[2]$.
For an integer $h$ from $i+1$ to $j$,
	we recursively define a subset $\mathcal{P}^{\alp_{i}\cdots\alp_{h-1}}$ as
	\envHLineCmDef
	{
		\mathcal{P}^{\alp_{i}\cdots\alp_{h-1}}
	}
	{	
		\mpN{ \SetT{ p |_{x_{h-1} = \alp_{h-1} } }{ p \in \mathcal{P}^{\alp_{i}\cdots\alp_{h-2}} } }	
	}	
	where $\mathcal{P}^{\alp_{i}\cdots\alp_{h-2}} = \mathcal{P}$ if $h=i+1$.
Then we have
	\envHLinePdPt[5_Lem2_Eq1]{\subset}
	{
		\mathcal{P}^{\alp_{i}\cdots\alp_{j-1}}
	}
	{
		\mpN{ \SetT{ p |_{x_{i} = \alp_{i} } |_{x_{i+1} = \alp_{i+1} } \cdots   |_{x_{j-1} = \alp_{j-1} }  }{ p \in \mathcal{P}  } }
	}
\end{lemma}

The proofs of the lemmas will be given after that of the proposition.

\envProof[\refProp{5_Prp1}]{	
It immediately follows 
	from \refEq{4_Lem1_Eq1} and \refEq{5_Prp1_FmlQ}
	that
	\envOTLineThCm
	{
		\nmSet{ \fmQ{j} }
	}
	{
		2^{j+1-I_j} - 2
	}
	{
		2^{\paBW+1}-2
	}
	which proves \prpS{1}.

We will show		
	\prpS{3}, \prpS{4} and \refEq{5_Prp1_Main2Eq}  
	by induction on $j$ from $\paBW+1$ to $n+1$.
Suppose that $j=\paBW+1$.
We do not need 
	to prove \prpS{4} 
	because 
	$\fmQ{\paBW+1}$ is set to $\fmP{\paBW+1}$ by the initial condition
	and
	calculation is unnecessary.
We can easily  
	verify \refEq{5_Prp1_Main2Eq}  
	because 
	it is equal to \refEq{3_Thm1_Main2Eq} in \refThm{3_Thm1}.
We will prove \prpS{3}.
We may assume $i=1$,
	since
	$I_{\paBW+1}=1$
	and
	\prpS{3} for $i>I_{\paBW+1}$ holds by \prp{3} in \refThm{3_Thm1}. 	
From \refEq{5_Lem1_Eq1} with $(i,j)=(1,\paBW+1)$ 
	and
	\refEq{5_Lem2_Eq1} with $\mathcal{P}=\fsP{1}{\tplE}$,	
	we see that
	\envPLinePd[5_Pr_Prp1_Eq2]
	{
		\fsP{\paBW+1}{\alp_{1}\cdots\alp_{\paBW}} 	
	\lnP{\subset} 	
		\mpN{ \rgBP_{\paBW+1,\paBW+1} }
	\lnP{=} 	
		\mpN{ \rgBP[x_{\paBW+1}] }
	}
We define
	\envHLine
	{
		\fmX_{x}
	}
	{ 
		\SetO{ \phi, \SetO{1}, \Set{x}, \Set{\lnot x}, \Set{x, \lnot x} }
	}	
	for a variable $x$.
By
	\refEq{3_PL_DefMapN},	
	any subset of $\mpN{ \rgBP[x_{\paBW+1}] }$ 
	belongs to $\fmX_{x_{\paBW+1}}$,
	which,
	together with \refEq{5_Pr_Prp1_Eq2} and $\fsQ{\paBW+1}{\alp_{1}\cdots\alp_{\paBW}} = \fsP{\paBW+1}{\alp_{1}\cdots\alp_{\paBW}}$,
	implies
	\envHLinePdPt[5_Pr_Prp1_Eq3]{\in}
	{
		\fsQ{\paBW+1}{\alp_{1}\cdots\alp_{\paBW}} 
	} 
	{
		\fmX_{x_{\paBW+1}}
	} 
Thus
	\envMCmPt{\leq}
	{
		\nmSet{ \fsQ{\paBW+1}{\alp_{1}\cdots\alp_{\paBW}} } 
	}
	{
		2
	}
	and
	we prove \prpS{3} for $i=1$.	

Suppose that 
	$j>\paBW+1$,
	and
	\prpS{3}, \prpS{4} and \refEq{5_Prp1_Main2Eq} are true in the case of $j-1$.
	
Firstly
	we will prove \prpS{3} for the case of $j$.
Let $\fmQ[m]{j}$ 
	denote 
	the family consisting of 
	the temporal subsets defined in \refEq{5_Prp1_DefQm}:	
	\envHLinePd[5_Pr_Prp1_FmlQm]
	{
		\fmQ[m]{j} 
	}
	{ 
		\SetT{  \fsQ[m]{j}{\alp_{i}\cdots\alp_{j-1}} }{ I_{j-1} \leq i \leq j-1, (\alp_{i}, \cdots, \alp_{j-1}) \in \setF[2]^{j-i} }
	}
The following properties hold.	
\envItem{
\item[\prpp{3}]
	\envMPdPt{\leq}
	{
		\nmSet{ \fsQ[m]{j}{\alp_{i}\cdots\alp_{j-1}} }		
	}
	{		
		\nmSet{ \fsQ{j-1}{\alp_{i}\cdots\alp_{j-2}} }
	}	
\item[\prpp{4}]
The computing time of \refEq{5_Prp1_DefQm} for all elements of $\fmQ[m]{j}$ 
	is 
	bounded by
	\envPLinePd
	{
		\sbL[a]{ l \SmT[l]{i=I_{j-1}}{j-1} \Sm[l]{ (\alp_{i}, \ldots, \alp_{j-2}) \in \setF[2]^{j-1-i} } \nmSet[a]{ \fsQ{j-1}{\alp_{i}\cdots\alp_{j-2}} }  }
	}
}
These properties can be 
	shown 
	as in the cases of \prp{3} and \prp{4} in \refThm{3_Thm1}.		
We omit their proofs for space limitation.\footnote{
We give brief explanations.
Both definitions of $\fsP{j}{\alp_{i}\cdots\alp_{j-1}}$ and $\fsQ[m]{j}{\alp_{i}\cdots\alp_{j-1}}$ are almost same
	as we see from \refEq{3_Thm1_DefP} and \refEq{5_Prp1_DefQm};
	only the conditions $p\in \fsP{j-1}{\alp_{i}\cdots\alp_{j-2}}$ and $p\in \fsQ{j-1}{\alp_{i}\cdots\alp_{j-2}}$ differ.
We also see 
	from \refEq{3_Thm1_FmlP} and \refEq{5_Pr_Prp1_FmlQm}
	that
	both definitions of $\fmP{j}$ and $\fmQ[m]{j}$ are almost same;	
	the conditions $1 \leq i \leq j-1$ and $I_{j-1} \leq i \leq j-1$ differ.
We can prove \prpp{3} and \prpp{4} 
	in the same ways as \prp{3} and \prp{4},
	respectively,	
	by commuting the above different places. 
}	
We will prove \prpS{3}.
Let $\fsQ{j}{\alp_{i}\cdots\alp_{j-1}} \in \fmQ{j}$.
If $i>I_j$,
	we see from \refEq{5_Prp1_DefQ}
	that
	\envHLineCm
	{
		\fsQ{j}{\alp_{i}\cdots\alp_{j-1}} 
	}
	{ 
		\fsQ[m]{j}{\alp_{i}\cdots\alp_{j-1}}
	}
	which, 
	together with \prpp{3},	
	yields
	\envHLinePdPt[5_Pr_Prp1_Eq5]{\leq}
	{
		\nmSet[a]{ \fsQ{j}{\alp_{i}\cdots\alp_{j-1}} }		
	}
	{
		\nmSet[a]{ \fsQ{j-1}{\alp_{i}\cdots\alp_{j-2}} }		
	} 
Assume $i = I_j$.
We also see from \refEq{5_Prp1_DefQ}
	that	
	\envPLinePd
	{
		\fsQ{j}{\alp_{I_j}\cdots\alp_{j-1}} 
	\lnP{=} 
		\SetO{1}
	\quad\text{or}\quad
		\fsQ[m]{j}{\alp_{I_j}\cdots\alp_{j-1}}
	}	
By \refEq{5_Lem1_Eq1} with $i = I_j$
	and
	\refEq{5_Lem2_Eq1} with $\mathcal{P}=\fsP{I_j}{\tplE}$,	
	we obtain
	\envMOCm
	{
		\fsQ[m]{j}{\alp_{I_j}\cdots\alp_{j-1}}  
	\lnP{\subset} 	
		\mpN{ \rgBP[x_{j}] }
	}
	and
	\envMPdPt{\in}
	{
		\fsQ[m]{j}{\alp_{I_j}\cdots\alp_{j-1}}  
	}
	{
		\fmX_{x_j}
	}	
Thus	
	$\fsQ{j}{\alp_{I_j}\cdots\alp_{j-1}} \in \fmX_{x_j}$,
	and
	\envHLinePdPt[5_Pr_Prp1_Eq6]{\leq}{
		\nmSet[a]{ \fsQ{j}{\alp_{I_j}\cdots\alp_{j-1}}  } 
	} 
	{
		2
	}
It follows from \refEq{5_Pr_Prp1_Eq5} and \refEq{5_Pr_Prp1_Eq6}
	that
	\prpS{3} holds in the case of $j$. 
	
Next
	we will prove \prpS{4} for the case of $j$.
When $i>I_j$,
	setting $\fsQ{j}{\alp_{i}\cdots\alp_{j}} = \fsQ[m]{j}{\alp_{i}\cdots\alp_{j}}$
	is only required
	in \refEq{5_Prp1_DefQ}.
Thus,
	the computing time of \refEq{5_Prp1_DefQm} and \refEq{5_Prp1_DefQ} 
	for the elements $\fsQ{j}{\alp_{i}\cdots\alp_{j-1}}$ with $i>I_j$
	is bounded by
	the time stated in \prpp{4}.
Therefore,
	to prove \prpS{4},
	we may show that
	the computing time of \refEq{5_Prp1_DefQ} 
	for all elements $\fsQ{j}{\alp_{i}\cdots\alp_{j-1}}$ with $i=I_j$
	is bounded by $\sbL{2^{\paBW}}$.
Let $(\alp_{I_j}, \ldots, \alp_{j-1}) \in \setF[2]^{j-I_j}$
	and		
	let $\alp = \alp_{I_{j -1}} \in \setF[2]$.		
By \refEq{5_Lem1_Eq1} with $i=I_{j -1}$
	and
	\refEq{5_Lem2_Eq1} with $\mathcal{P}=\fsP{I_{j -1}}{\tplE}$,	
	we obtain
	\envMOCm
	{
		\fsQ{j-1}{\alp \alp_{I_{j}}\cdots\alp_{j-2}} 
	\lnP{\subset}
		\mpN{ \rgBP[x_{j-1}] }
	}	
	and 
	\envHLineCmPt[5_Pr_Prp1_Eq8]{ \osTx{\subset}{\refEq{5_Prp1_DefQm}}}
	{	
		\fsQ[m]{j}{\alp \alp_{I_{j}}\cdots\alp_{j-1}} 
	}
	{ 
		\mpN{ \SetT{ p |_{x_{j-1} = \alp_{j-1} } }{ p \in \rgBP[x_{j-1}] } }
	\lnP{\osTx{=}{\refEq{3_PL_EqBooPolSub}}}
		\mpN{ \setF[2] }
	\lnP{\osTx{\subset}{\refEq{3_PL_DefMapN}}}
		\SetO{1}
	} 
	which implies that
	$\fsQ[m]{j}{\alp \alp_{I_{j}}\cdots\alp_{j-1}}$ is either $\phi$ or $\SetO{1}$,
	or equivalently, 
	$\nmSet{\fsQ[m]{j}{\alp \alp_{I_{j}}\cdots\alp_{j-1}}}$ is either zero or one.
Hence,	
	the time of checking 
	whether 
	$\Pd[l]{\alp \in \setF[2]} \nmSet[a]{\fsQ[m]{j}{\alp \alp_{I_{j}}\cdots\alp_{j-1}} }$
	is zero or not is in $\sbL{1}$.
By the definition of \refEq{5_Prp1_DefQ},
	we can calculate	
	the single element $\fsQ{j}{\alp_{I_j}\cdots\alp_{j}}$	
	in time $\sbL{1}$.
Since $\nmSet{\setF[2]^{j-I_j}} = \nmSet{\setF[2]^{\paBW}} = 2^{\paBW}$,		
	the computing  time of \refEq{5_Prp1_DefQ} 
	for all elements $\fsQ{j}{\alp_{i}\cdots\alp_{j-1}}$ with $i=I_j$	
	is bounded by $\sbL{2^{\paBW}}$.	
This concludes that
	\prpS{4} in the case of $j$ is true.

Finally 
	we will prove \refEq{5_Prp1_Main2Eq} for the case of $j$.
By the induction hypothesis, 
	$\plP[]{j-1}{}$ satisfies \refEq{5_Prp1_Main2Eq}.
Using the distributivity of $\lor$ over $\land$,	 
	we can obtain
	the following equation 	
	as in \refEq{4_Pr_3_Thm1_Eq2} and \refEq{4_Pr_3_Thm1_Eq3}:	
	\envHLinePd[5_Pr_Prp1_FML_Eq1]
	{
		\plP[]{j}{} 
	}
	{
		\bkR[g]{ \Ad{\alp_{j-1} \in \setF[2]} 
			\bkR[g]{\cdots
				\bkR[g]{
					\Ad{\alp_{I_{j}} \in \setF[2]}
					\bkR[g]{
						\Ad{\alp_{I_{j-1}} \in \setF[2]} \plF[]{\fsQ[m]{j}{\alp_{I_{j-1}} \cdots \alp_{j-1}} }{} 
					}\lor \plF[]{\fsQ[m]{j}{\alp_{I_{j}} \cdots \alp_{j-1}}}{} 
				}
			\cdots  } \lor \plF[]{ \fsQ[m]{j}{\alp_{j-1}} }{}
		} 
		\lnAH
		\lor \plF[]{ \fsP{j}{\tplE} \cup \cdots \cup \fsP{n+1}{\tplE} }{} 		
	}
Because	
	\envMCmPt{\in}
	{
		\fsQ[m]{j}{\alp \alp_{I_{j}}\cdots\alp_{j-1}} 
	}
	{
		\SetO{ \phi, \SetO{1} }
	}	
	we see from \refEq{2_PL_DefFacFrm} 
	that	
	\envHLineCm
	{
		\plF[]{\fsQ[m]{j}{\alp\alp_{I_{j}} \cdots \alp_{j-1}} }{}  	
	}
	{
		\nmSet[a]{\fsQ[m]{j}{\alp \alp_{I_{j}}\cdots\alp_{j-1}} } 	
	}
	where 
	the values $0$ and $1$ in $\setF[2]$ are identified with those in $\setZ$.
Hence
	\envHLineCm
	{
		\AD{\alp_{I_{j-1}} \in \setF[2]} \plF[]{\fsQ[m]{j}{\alp_{I_{j-1}} \alp_{I_{j}} \cdots \alp_{j-1}} }{} 		
	}
	{
		\Pd{\alp_{I_{j-1}} \in \setF[2]} \nmSet[a]{\fsQ[m]{j}{\alp_{I_{j-1} } \alp_{I_{j}}\cdots\alp_{j-1}} }	
	}	
	and
	\envLine
	{
		\bkR[g]{
			\Ad{\alp_{I_{j-1}} \in \setF[2]} \plF[]{\fsQ[m]{j}{\alp_{I_{j-1}} \alp_{I_{j}} \cdots \alp_{j-1}} }{} 
		}\lor \plF[]{\fsQ[m]{j}{\alp_{I_{j}} \cdots \alp_{j-1}}}{} 
	}
	{
		\envCaseTPd{
			1 \lor \plF[]{\fsQ[m]{j}{\alp_{I_{j}} \cdots \alp_{j-1}}}{}
			&
			\text{if $\Pd[l]{\alp_{I_{j-1}} \in \setF[2]} \nmSet[a]{\fsQ[m]{j}{\alp_{I_{j-1} } \alp_{I_{j}}\cdots\alp_{j-1}} } = 1$}	
		}{
			0 \lor \plF[]{\fsQ[m]{j}{\alp_{I_{j}} \cdots \alp_{j-1}}}{}
			&
			\text{if $\Pd[l]{\alp_{I_{j-1}} \in \setF[2]} \nmSet[a]{\fsQ[m]{j}{\alp_{I_{j-1} } \alp_{I_{j}}\cdots\alp_{j-1}} } = 0$}
		}
	}
By the annihilator and identity laws for $\lor$,
	\envLineThPd[5_Pr_Prp1_FML_Eq3]
	{
		\bkR[g]{
			\Ad{\alp_{I_{j-1}} \in \setF[2]} \plF[]{\fsQ[m]{j}{\alp_{I_{j-1}} \alp_{I_{j}} \cdots \alp_{j-1}} }{} 
		}\lor \plF[]{\fsQ[m]{j}{\alp_{I_{j}} \cdots \alp_{j-1}}}{} 
	}
	{
		\envCaseTCm{
			1 
			&
			\text{if $\Pd[l]{\alp_{I_{j-1}} \in \setF[2]} \nmSet[a]{\fsQ[m]{j}{\alp_{I_{j-1} } \alp_{I_{j}}\cdots\alp_{j-1}} } = 1$}	
		}{
			\plF[]{\fsQ[m]{j}{\alp_{I_{j}} \cdots \alp_{j-1}}}{}
			&
			\text{if $\Pd[l]{\alp_{I_{j-1}} \in \setF[2]} \nmSet[a]{\fsQ[m]{j}{\alp_{I_{j-1} } \alp_{I_{j}}\cdots\alp_{j-1}} } = 0$}
		}
	}
	{
		\plF[]{\fsQ{j}{\alp_{I_{j}} \cdots \alp_{j-1}} }{} 
	}	
Since $\fsQ[m]{j}{\alp_{i}\cdots\alp_{j}} = \fsQ{j}{\alp_{i}\cdots\alp_{j}}$ for $i>I_j$,
	combining 
	\refEq{5_Pr_Prp1_FML_Eq1} and \refEq{5_Pr_Prp1_FML_Eq3}
	gives
	\refEq{5_Prp1_Main2Eq} in the case of $j$.
	
We conclude 
	that
	all of \prpS{3}, \prpS{4} and \refEq{5_Prp1_Main2Eq} are true in the case of $j$,
	and
	we complete the induction step.
Therefore 
	\refProp{5_Prp1} holds.
}

\envProof[\refLem{5_Lem1}]{
The definition of $i$
	implies
	$p \in \rgBP_{i, n}$,
	and
	that of $\paBW$ 
	implies
	$p \in \rgBP_{i, i+\paBW}$.
Therefore,
	by \refEq{3_PL_EqBooPolSub},
	we obtain	
	\refEq{5_Lem1_Eq1}.
}
\envProof[\refLem{5_Lem2}]{
Obviously,
	\refEq{5_Lem2_Eq1} with $j=i+1$ holds by definition.
Since 
	$\mpN{\mathcal{Q}} \subset \mathcal{Q}$
	for any subset $\mathcal{Q}$ of $\rgBP_n$,
	\refEq{5_Lem2_Eq1} with $j=i+2$ 
	is 
	proved by
	\envHLineFPd
	{
		\mathcal{P}^{\alp_{i}\alp_{i+1}}
	}
	{
		\mpN{ \SetT{ p |_{x_{i+1} = \alp_{i+1} } }{ p \in \mathcal{P}^{\alp_{i}} } }
	}
	{
		\mpN{ \SetT{ p |_{x_{i+1} = \alp_{i+1} } }{ p \in \mpN{ \SetT{ p |_{x_{i} = \alp_{i} } }{ p \in \mathcal{P} } } } }
	\lnAHP{\subset}
		\mpN{ \SetT{ p |_{x_{i+1} = \alp_{i+1} } }{ p \in \SetT{ p |_{x_{i} = \alp_{i} } }{ p \in \mathcal{P} } } }
	}
	{
		\mpN{ \SetT{ p |_{x_{i} = \alp_{i} } |_{x_{i+1} = \alp_{i+1} } }{ p \in \mathcal{P} } }
	}
Similarly,
	we can prove \refEq{5_Lem2_Eq1} for general $j$
	using induction on $k=j-i$.
We omit the details for space limitation.
}

\section{Proofs of \refCor[s]{3_Cor1} and \ref{3_Cor2}}	\label{sectSix}
We will prove \refCor{3_Cor1}.
\envProof[\refCor{3_Cor1}]{	
Suppose that $\syN$ is BMQ.
Then
	the Boolean polynomials in $\syN$ are quadratic polynomials,
	and
	their degrees are at most $2$.
The number of monomials of degrees at most $2$ is bounded by $\sbL{n^2}$,
	and 
	$\syN$ belongs to $\clST{n}{l}$ with $l=\sbL{n^2}$.
Thus \refThm{3_Thm2} 
	implies 
	{\bf (i)}.

Suppose that $\syN$ is CNF.
Then
	the Boolean polynomials in $\syN$ are clause polynomials,
	and
	$\syN$ belongs to $\clST{n}{l}$ with $l=1$,
	which,
	together with \refThm{3_Thm2}, 
	proves {\bf (ii)}.
}

We will introduce some notions and facts for graph list-coloring 
	to prove \refCor{3_Cor2}.

Let $G$ be a graph,
	and
	let $\bkS{l}=\SetO{0, \ldots, l-1}$ be $l$ distinct colors.
We suppose $l \geq 2$
	because
	list coloring problem of one color is trivial.
We mean by $k$ the integer such that
	$2^{k-1} < l \leq 2^k$,
	and
	consider $\bkS{2^k}$ a universal set of colors.	
For a $\bkS{l}$-list $L$,
	we denote 
	by $\dlC{L}$
	the complement of $L$,
	i.e.,
	$\dlC{L} = \bkS{2^k} \setminus L$.
We define
	a bijection from $\setF[2]^k$ to $\bkS{2^k}$ by
	\envHLineCmDef
	{
		\mpF{\gam_1, \ldots, \gam_k}
	}
	{
		\gam_1  + \gam_2 2 +  \cdots + \gam_k 2^{k-1} 
	}
	where
	$\gam_1, \ldots, \gam_k$ are $k$ values in $\setF[2]$.
Let 
	$x_{v,1}, \ldots, x_{v,k}$
	be
	$k$ variables 
	associated with a vertex $v$.
For a color $c$ in $\bkS{2^k}$,
	we define 
	a clause polynomial of degree $k$ by
	\envHLineCmDef
	{
		\clG[]{v,c}{}	
	}
	{
		\PdT{h=1}{k} \bkR{ x_{v,h} + \gam_h + 1} 
	\lnP{\in}
		\rgBP[x_{v,1}, \ldots, x_{v,k}]
	}
	where
	$(\gam_1, \ldots, \gam_k) = \mpFi{c}$.
For an edge $uv$,
	we also define
	\envHLineDefPd
	{
		\clG[]{uv,c}{}	
	}
	{
		\PdT{h=1}{k} \bkR{ x_{u,h} + x_{v,h} + \gam_h + 1} 
		\lnP{\in}
		\rgBP[x_{u,1}, \ldots, x_{u,k},x_{v,1}, \ldots, x_{v,k}]
	}
Let $\syN_{(G,L)}$  
	be
	a system in the variables $x_{v,h}$	
	which consists of 
	the following equations:
	\envPLineCm[]
	{\label{6_PL_defSysV}
		\clG[]{v, c}{} &\lnP{=}& 0	
		\qquad
		(v \in V, c \in \dlC{L_v})
	,\\\label{6_PL_defSysE}
		\clG[]{uv, 0}{} &\lnP{=}& 0	
		\qquad
		(uv \in E)
	}	
	where
	$L_v$ are $\bkS{l}$-lists of allowed colors for vertices $v$.
Note that
	the color $0$ 		
	is corresponding to the zero tuple,		
	and
	\envMPd
	{
		\clG[]{uv,0}{}	
	}
	{
		\PdT{h=1}{k} \bkR{ x_{u,h} + x_{v,h} + 1} 
	}
	
Let 
	$\bfAlp_v \in \setF[2]^k$
	and
	set
	$a_v = \mpF{\bfAlp_v}$
	for vertices $v$.
We see that
	$a_v \neq c$
	if and only if 
	$\clG{v, c}{\bfAlp_v} = 0$ for a color $c$,
	because
	zero 
	is
	an annihilating element for product.
We also see that
	$\clG{uv, 0}{\bfAlp_u,\bfAlp_v} = \clG{u, a_v}{\bfAlp_u} = \clG{v, a_u}{\bfAlp_v} $
	because of the definitions.
By these facts
	we can find the following properties:
\envItem{
\item[(V)]
For a vertex $v$,
	the color $a_v$ is in $L_v$ 
	if and only if
	$\clG{v, c}{\bfAlp_v} = 0$ for all $c \in \dlC{L_v}$.
\item[(E)]
For an edge $uv$,
	the colors $a_u$ and $a_v$ are different 
	if and only if
	$\clG{uv, 0}{\bfAlp_u,\bfAlp_v} = 0$.
}

We will show that
	the list-colorability of $G$ on $(L_v)_{v \in V}$
	is 
	equivalent to
	the satisfiability of $\syN_{(G,L)}$.
Suppose that
	$G$ is list-colorable.
Then 
	there exists a tuple $(a_v)_{v\in V}$ of colors 
	such that
	(i) $a_v \in L_v$ for every $v \in V$;
	and
	(ii) $a_u \neq a_v$ for every $uv \in E$.
It follows 
	from (i) and (V)
	that
	$(\mpFi{a_v})_{v\in V}$ satisfies \refEq{6_PL_defSysV},
	and
	from (ii) and (E)
	that
	$(\mpFi{a_v})_{v\in V}$ satisfies \refEq{6_PL_defSysE}.
Hence 
	$(\mpFi{a_v})_{v\in V}$ is a solution,
	and
	$\syN_{(G,L)}$ is satisfiable.
Suppose that
	$\syN_{(G,L)}$ is satisfiable,
	and
	$(\bfAlp_v)_{v \in V}$ is its solution.
Similarly to the above,
	it can be seen that
	$(\mpF{\bfAlp_v})_{v \in V}$ is a proper assignment of colors.
Thus
	$G$ is list-colorable. 

We are in a position to prove \refCor{3_Cor2}.

\envProof[\refCor{3_Cor2}]{
We may assume that
	$l \geq 2$
	and 
	$k \geq 1$,
	where
	$2^{k-1} < l \leq 2^k$.
Let $\syN_{(G,L)}$  
	be
	the system 	
	defined by \refEq{6_PL_defSysV} and \refEq{6_PL_defSysE}.	
The Boolean polynomials in the system  
	are in
	$\rgBP[ (x_{v,h})_{v \in V, 1 \leq h \leq k}]$,
	and
	$\syN_{(G,L)} \in \clST{kn}{2^{kn}}$.
	
Firstly,
	we will show 
	\envHLinePdPt[6_Pr_3_Cor2_Eq1]{\in}
	{
	 	\syN_{(G,L)} 
	}
	{
	 	\clST{kn}{2l}
	}
If $g$ is $\clG[]{v,c}{}$ in \refEq{6_PL_defSysV},
	then
	$g$ is a clause and 
	$g \in \spCLT{kn}{1} \subset \spCLT{kn}{2l}$.
Suppose that  $g$ is $\clG[]{uv,0}{}$ in \refEq{6_PL_defSysE}.		
Then
	\envHLineThPd
	{
		g	
	}
	{
		\PdT{h=1}{k} \bkR{ x_{u,h} + x_{v,h} + 1} 
	}
	{
		\tpSm{ H_u,  H_v \subset \Set{1, \ldots, k}}{ \nmSet{H_u} + \nmSet{H_v} = k, H_u \cap H_v = \phi} 
		\bkR[a]{ \Pd{ h_u \in H_u } x_{u,h_u} }
		\bkR[a]{ \Pd{ h_v \in H_v }\bkR{ x_{v,h_v}  + 1 } }
	\lnP{\in}
		\spCLT{kn}{2^k}
	}
Because $2^{k-1} < l$,
	we have
	$2^k < 2l$
	and
	$\spCLT{kn}{2^k} \subset \spCLT{kn}{2l}$.
Therefore
	$g \in \spCLT{kn}{2l}$.
Since
	$\syN_{(G,L)}$ consists of Boolean polynomials in \refEq{6_PL_defSysV} and \refEq{6_PL_defSysE},
	we obtain \refEq{6_Pr_3_Cor2_Eq1}.

Let $v_1, \ldots, v_n$
	be 
	vertices whose order give the bandwidth $\paBW$.
Referring to the order of vertices,
	we define  
	that of variables by 	
	\envPLineCm
	{
		x_{1,1}, \ldots, x_{1,k},
		x_{2,1}, \ldots, x_{2,k},
		\ldots,
		x_{n,1}, \ldots, x_{n,k}
	}
	where
	$x_{i,h} = x_{v_i,h}$.
By \refEq{3_PL_DefWthP} and \refEq{3_PL_DefWthG},
	the bandwidth of $\syN_{(G,L)}$
	is 
	$(\paBW+1) k -1$.
With \refEq{6_Pr_3_Cor2_Eq1},
	\refThm{3_Thm2}
	implies that
	the satisfiability of $\syN_{(G,L)}$
	is 
	decidable in time
	\envPLineCm
	{
		\sbL{ 2^{(\paBW+1) k -1} (2l) (m' + kn) }
	}
	where
	$m'$
	is 
	the number of equations in \refEq{6_PL_defSysV} and \refEq{6_PL_defSysE}.
We have
	\envPLinePd
	{
		2^{(\paBW+1) k }	\lnP{<}	(2l)^{\paBW+1}
	,\qquad
		m'		\lnP{\leq}		2^k n + m 		\lnP{<}		2l n + m
	,\qquad
		kn		\lnP{<}		ln
	}
Thus 
	the time of solving $\syN_{(G,L)}$ is bounded by
	\envMOPd
	{
		\sbL{ (2l)^{\paBW+2} (m+ln) }
	}
This completes the proof,
	because
	the list-colorability of $G$ on $(L_v)_{v \in V}$
	is 
	equivalent to
	the satisfiability of $\syN_{(G,L)}$.
}


\section*{Acknowledgement}		
The author would like to thank Tomohiro Sonobe for his helpful comments.

\renewcommand{\thesubsubsection}{Appendix \Alph{subsubsection}}




\end{document}